\documentclass[11pt,twoside]{article}

\usepackage{epsfig,amsfonts,color}
\usepackage{amsmath, bbm, mathabx}
\usepackage[normalem]{ulem}

\usepackage{amssymb, palatino, geometry,url}
\usepackage{algorithmic}
\usepackage[noresetcount,lined,boxed]{algorithm2e} 

\usepackage[colorlinks=true,linkcolor=blue,citecolor=blue,urlcolor=blue]{hyperref}
\usepackage{subcaption}
\usepackage{amsthm}
\usepackage[utf8x]{inputenc}
\usepackage{textgreek}
\usepackage{multirow}

\usepackage{fancyhdr}
\usepackage{lineno}
\usepackage{float}
\usepackage[section]{placeins}
\usepackage[table]{xcolor}

\title{SAFE-OCC: A Novelty Detection Framework for Convolutional Neural Network Sensors and \\its Application in Process Control}
\author{Joshua L. Pulsipher$^1$\;, Luke D. J. Coutinho$^1$, Tyler A. Soderstrom$^2$, and Victor M. Zavala$^1$\thanks{Corresponding Author: victor.zavala@wisc.edu}\\
	{\small $^1$Department of Chemical and Biological Engineering}\\
	{\small \;University of Wisconsin-Madison, 1415 Engineering Dr, Madison, WI 53706, USA}\\
	{\small $^2$Advanced Process Control}\\
	{\small \;ExxonMobil Research and Engineering, 22777 Springwoods Village Pkwy, Spring, TX 77389, USA}}
\date{}

\begin{document}
	
\maketitle

\begin{abstract}
We present a novelty detection framework for Convolutional Neural Network (CNN) sensors that we call Sensor-Activated Feature Extraction One-Class Classification (SAFE-OCC). We show that this framework enables the safe use of computer vision sensors in process control architectures. Emergent control applications use CNN models to map visual data to a state signal that can be interpreted by the controller. Incorporating such sensors introduces a significant system operation vulnerability because CNN sensors can exhibit high prediction errors when exposed to novel (abnormal) visual data. Unfortunately, identifying such novelties in real-time is nontrivial. To address this issue, the SAFE-OCC framework leverages the convolutional blocks of the CNN to create an effective feature space to conduct novelty detection using a desired one-class classification technique. This approach engenders a feature space that directly corresponds to that used by the CNN sensor and avoids the need to derive an independent latent space. We demonstrate the effectiveness of SAFE-OCC via simulated control environments. 
\end{abstract}


\section{Introduction} \label{sec:intro}

Advances in photogrammetry technologies have led to a significant rise in the application of computer vision strategies \cite{hajirassouliha2018suitability}. Such implementations leverage computer vision methods to automatically extract desired information from image/visual data. This powerful capability has proved valuable in diverse fields such as autonomous vehicle environment sensing \cite{janai2020computer}, facial recognition \cite{balaban2015deep}, motion analysis \cite{poppe2007vision}, automatic inspection \cite{neethu2015role}, robotic navigation \cite{bonin2008visual}, biomedical image analysis \cite{suzuki2017overview}, and quality control \cite{wu2013colour,villalba2019deep}. However, these techniques have seen limited applications in process control, where it is common practice for operators to manually manipulate parts of a process based on visual information (i.e., form a manual control loop using human-perceived visual insights) \cite{villalba2019deep}. It is thus desired to leverage computer vision approaches to remove the human sensory component and with this design  automated control loops that incorporate sensors that directly extract state information from raw visual data. Such an approach is exemplified in recent work \cite{lu2020image}, where the authors  propose to incorporate image data into a model predictive control framework. The use of computer vision sensors  can lead to enhanced safety, increased efficiency, and reduced operational costs \cite{antsaklis1991introduction}.  
\\

Computer vision is often employed in the form of supervised learning approaches which entail two principal categories: classification and regression \cite{hastie2009overview}. Classification schemes seek to map an input image to a descriptive category (e.g., object detection), whereas regression schemes map input images to real-valued (continuous) descriptors (e.g., state estimation). These denote supervised learning strategies since the predictor models are trained via (learn from) labeled image data. Training data typically contain a large set of representative image-output pairs where the labeled outputs denote the desired information (i.e., states) that should be extracted from each image. Image augmentation is typically employed to increase the size and generality of the training set by means of augmenting images with synthetic perturbations (disturbances) such as rotation, stretching, translation, splattering, noise, and more \cite{shorten2019survey}. The augmented set images help to generalize the prediction domain of supervised learning models and avoids rotation variance issues. Once trained, these predictors can rapidly estimate the states of new images under the assumption that they are members of the prediction domain spanned by the training set \cite{jiang2021convolutional}. 
\\

Convolutional Neural Networks (CNNs) are one of the most prevalent supervised learning models used in computer vision to extract information from visual data \cite{khan2018guide}. Their prevalence can be largely attributed to the ability of CNNs to automate feature extraction, eliminating the need for human-engineered filters (i.e., feature patterns) \cite{gu2018recent}. These machine learning models have also demonstrated their ability to match or exceed human performance with regular images (those spanned by the training set), and to facilitate much higher data throughput \cite{geirhos2017comparing}. CNNs provide new and exciting functionalities that can enable powerful process control applications and can also be used for diverse applications of interest to the process systems community (that go beyond computer vision), as highlighted recently in \cite{jiang2021convolutional}. Moreover, CNNs can nowadays be easily incorporated in computational workflows, since there are several open-source tools for implementing these models such as \texttt{Keras} \cite{geron2019hands}, \texttt{PyTorch} \cite{paszke2019pytorch}, and \texttt{Flux} \cite{innes2018fashionable}. 
\\

CNN models can be broken down conceptually into a couple of fundamental blocks: a feature extractor and a predictor \cite{khan2018guide}. An input image   is fed to the feature extractor which uses specialized convolutional layers that leverage filters (kernels) to extract informative visual patterns that derive meaningful predictive features. These extracted features are then given to a feedforward (dense) neural network predictor that predicts the desired output states. Model training involves simultaneously learning pattern filters for feature extraction and dense layer weights for the neural network predictor that minimize a loss function. 
\\

Although CNNs can predict to high accuracy, they typically incur large prediction errors on abnormal images (e.g., visual disturbances not accounted for in the training data) \cite{lee2017training}. This behavior is consistent with that of other learning models that typically excel at interpolation within the prediction domain, but can exhibit poor performance with input data that extrapolate beyond the span of the training data \cite{kosanovich1996improving}. This extrapolation data is often referred to as novel data, anomalies, or out-of-distribution data (here we use the term {\em novel}). Distinguishing whether an image set is novel or not (relative to the training data) is nontrivial in general since the data is high-dimensional in nature, training sets are usually large, and transitions into novelty may progress slowly over time. This introduces significant complications in incorporating computer vision sensors into real-time control loops, where it is vital that sensors provide accurate measurements in order to preserve robustness, stability, and efficiency. For this reason, traditional process control sensors (e.g., thermocouples and flow meters) are typically engineered to return an error signal when they malfunction (such that operators can take appropriate recourse action). It is thus critical that CNN sensors integrated within control systems are able to accurately identify novel images in real-time and return an error signal when appropriate. In other words, we require a real-time autonomous system to accurately assess the quality of the predictions made by a CNN sensor. 
\\

Novelty detection denotes a set of unsupervised learning methods that differentiate between novel and normal data. These methods denote an active area of research; Ruff and co-workers recently provided a thorough review in \cite{ruff2021unifying}. A couple of particularly prevalent paradigms are reconstruction models and one-class classification. Reconstruction model approaches seek to learn a low-dimensional latent space for a set of unlabeled training data that is predominately normal. New data is then encoded into the learned latent space, decoded back to its original dimension, and the reconstruction error incurred by this transformation is used as a metric for novelty (i.e., low reconstruction error for normal data and high reconstruction error for novel data). In the context of computer vision, this is typically done with encoders and decoders learned from convolutional autoencoder and/or variational autoencoder models. Some recent works include the probabilistic reconstruction approach using variational autoencoders in \cite{an2015variational}, the autoecoder ensemble approach in \cite{chen2017outlier}, the generative adversarial network based approach in \cite{zenati2018efficient}, and the adversarial mirrored autoencoder approach in \cite{somepalli2020unsupervised}. This class of methods provides rich capabilities for novelty detection with general datasets; however, in our context,  their principal disadvantage is that their latent space is independent of the feature space used by the CNN sensor. As such, these methods may be overly conservative in identifying novelties that may have limited effect on the CNN  performance. Moreover, they tend to be computationally expensive/complex to implement since they typically require the training of additional deep learning models (e.g., convolutional autoencoders). 
\\

One-Class Classification (OCC) denotes an area of methods that learn a single class of normal instances from unlabeled training data (typically assumed to be comprised of normal instances). These then identify novel data instances by determining if they lie outside the learned class. Here, a large focus of recent work has been on deep end-to-end approaches that extract image features and learn the normal class. Such work includes the deep one-class classification presented in \cite{perera2019learning}, the $E^3Outlier$ approach presented in \cite{wang2019effective}, and the deep end-to-end one-class classifier presented in \cite{sabokrou2020deep}. However, like the reconstruction techniques, these use a latent space that is derived independently of the feature space employed by the CNN predictor, and they tend to be complex and expensive to train/setup. 
\\

One-Class Support Vector Machines (OC-SVMs) are another prevalent group of techniques employed and work by adapting support vector machines to learn a boundary around the training data that thresholds the distinction between novel and normal data. In particular, Support Vector Data Description (SVDD) is a common choice that learns a spherical boundary around the training data in its feature space. In contrast to deep learning methods, OC-SVMs require a minor computational expense to train and are easy to implement with tools like \texttt{SciKit-Learn} \cite{pedregosa2011scikit}. However, OC-SVMs are only amenable for one-dimensional feature data (i.e., each data point is a vector of features) and cannot be directly used on image data. A number of methods overcome this challenge by using deep learning techniques to extract features from image data that can be used with OC-SVMs. For instance, Erfani et. al. extract features by training deep belief networks \cite{erfani2016high}, and Andrews and co-workers use pre-trained CNNs to extract features \cite{andrews2016transfer}. These approaches are attractive with their general simplicity and low computational burden, but have not yet been explored with feature spaces that are derived from a target CNN sensor of interest for the purpose of assessing the quality of its predictions.

\begin{figure}[!htb]
    \centering
    \includegraphics[width=0.95\textwidth]{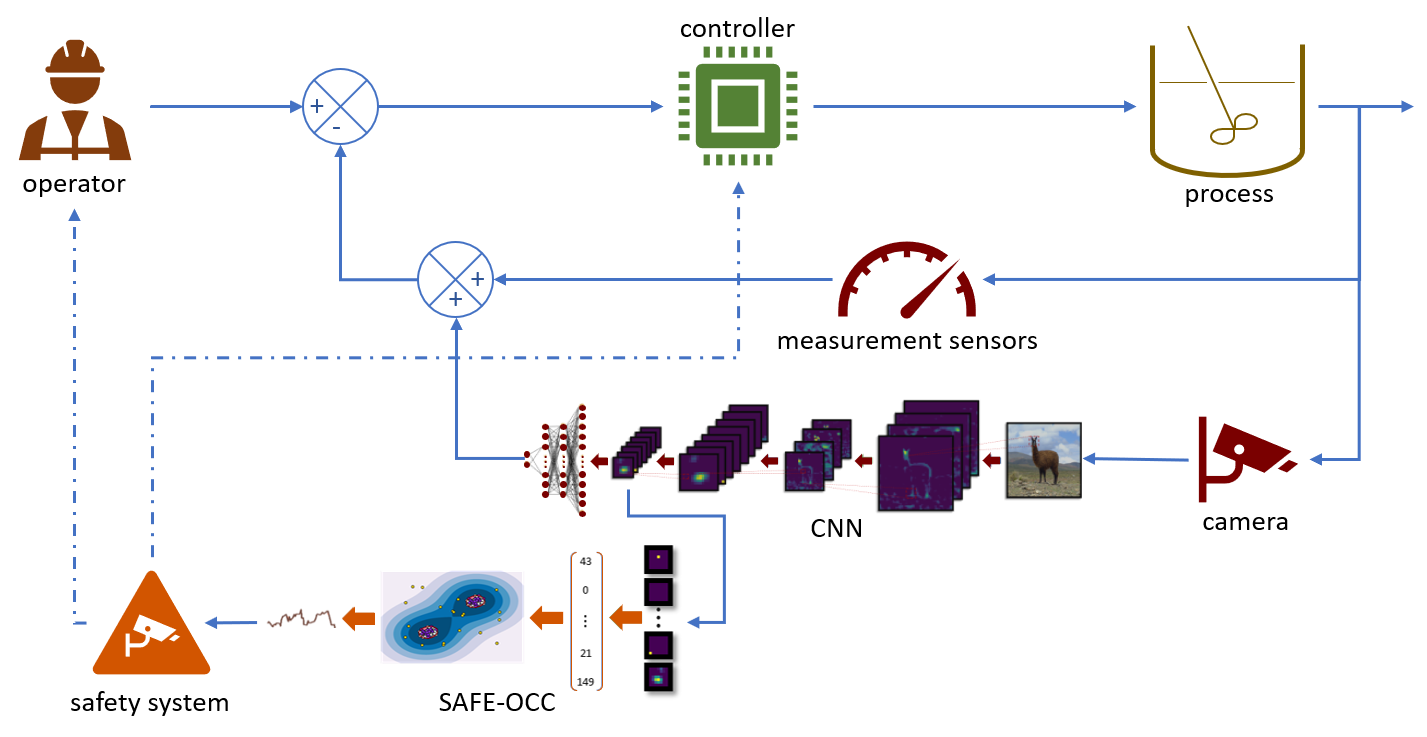}
    \caption{A visual abstract for the SAFE-OCC framework and its incorporation into process control architecture.}
    \label{fig:abstract}
\end{figure}

To address these limitations, here we propose the Sensor Activated Feature Extraction Once-Class Classification (SAFE-OCC) novelty detection framework. This framework assesses the novelty of input images relative to a trained CNN sensor by using a feature space directly derived from the CNN to conduct OCC as summarized in Figure \ref{fig:abstract}. This approach takes inspiration from and provides a more rigorous treatment of the conceptualized feature-based novelty detection approach outlined in \cite{exxon2020patent}. Specifically, we reduce targeted feature map outputs from the convolutional layers to produce a one-dimensional feature space which is amenable to standard OCC techniques (e.g., OC-SVM) to determine the novelty of input images. Our approach eliminates the need to create/train a feature extractor and instead derives a feature space that directly corresponds to the CNN sensor (making it less conservative in identifying novel images). These benefits make SAFE-OCC readily suited for a wide breadth of application fields that need to assess the prediction confidence made by a CNN sensor in real-time (including our applications of interest in process control).
\\

The paper is structured as follows. Section \ref{sec:background} provides relevant background and establishes notation for our framework. Section \ref{sec:compvis_control} details the incorporation of computer vision sensors in closed-loop control. Section \ref{sec:framework} details the SAFE-OCC framework. Section \ref{sec:examples} provides illustrative case studies using simulated control environments. Section \ref{sec:conclusion} details the key takeaways and plans for future work.

\section{Relevant Background and Notation} \label{sec:background}

In this section, we highlight pertinent background information and define notation for CNN models, PCA reduction techniques, OCC approaches, and process control. For a thorough review on these topics, we refer the reader to \cite{jiang2021convolutional}, \cite{zhang20052d}, \cite{khan2014one}, and \cite{seborg2010process}. 

\subsection{Convolutional Neural Networks} \label{sec:cnn}
CNNs denote a broad class of machine learning models that utilize specialized convolutional blocks to extract features from grid data objects. In this work, we will focus on CNN sensors $f_{\text{cnn}} : \mathbb{R}^{n_{v} \times n_v \times p} \mapsto \mathbb{R}^{n_y}$ that map an input image $V \in \mathbb{R}^{n_v \times n_v \times p}$ to state predictions $\hat{y} \in \mathbb{R}^{n_y}$. For simplicity in presentation, we consider square images $V$ with $p \in \mathbb{Z}_+$ color channels and indexing $V_{x_1,x_2,i}$ with $x_1,x_2 \in \{-N_v,N_v\}$ and $i \in \{1,p\}$, where the notation $\{a, b\}$ denotes the set of integers $a,\dots,b$ and we have $n_v = 2N_v + 1$; however, more general images and/or video can be used following the same techniques described in this work. 

\begin{figure}[!htb]
    \centering
    \includegraphics[width=0.8\textwidth]{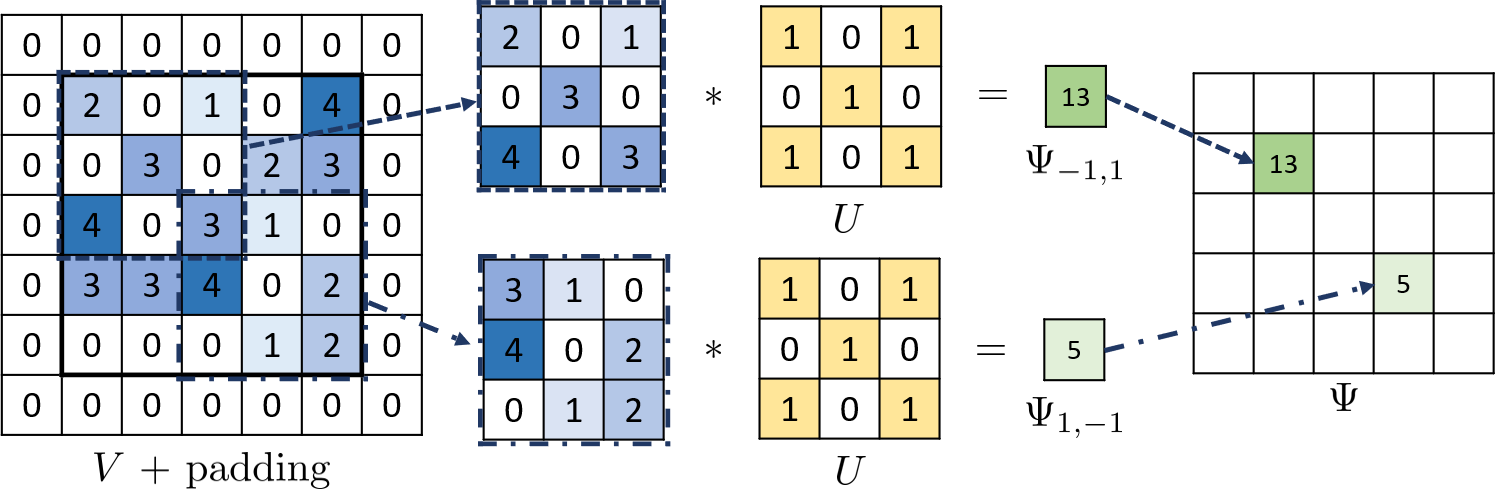}
    \caption{Convolution of a one-channel image $V$ with a single filter convolutional operator $U$ to produce feature map $\Psi$. Note that zero-padding is added to $V$.}
    \label{fig:convolution}
\end{figure}

Here a particular convolutional layer takes an input array $V$ and \emph{convolves} it by applying a convolutional operator $U \in \mathbb{R}^{n_u \times n_u \times p \times q}$ which employs $q \in \mathbb{Z}_+$ convolutional filters (i.e., kernels) to produce a feature map $\Psi \in \mathbb{R}^{n_\psi \times n_\psi \times q}$ (typically $n_\psi = n_v$). Again for simplicity in presentation, we consider square filters of size $n_u \times n_u \times p$, but this condition is readily relaxed. This operation is formalized as $\Psi = U * V$:
\begin{equation}
    \Psi_{x_1,x_2,j} = \sum_{i = 1}^p \sum_{x_1' = -N_u}^{N_u} \sum_{x_2' = -N_u}^{N_u}  U_{x_1',x_2',i,k} \cdot V_{x_1 + x_1', x_2 + x_2', i}
    \label{eq:convolution}
\end{equation}
for $x_1,x_2 \in \{-N_v,N_v\}$ and $j \in \{1,q\}$. We observe that this operation is fully decoupled over each filter (indexed by $j$). This operation is not well defined by the boundaries of $V$ as some indexing will violate the image domain, but this is typically resolved by adding padding of zero-valued entries around the image (a technique called zero-padding). We can also express this operation in the compact functional form $\Psi = f_c(V; U)$ where $f_c : \mathbb{R}^{n_v \times n_v \times p} \mapsto \mathbb{R}^{n_\psi \times n_\psi \times q}$. Moreover, CNN convolutional operators typically exhibit the property that $n_u \ll n_v$. Figure \ref{fig:convolution} illustrates the convolution of a $5 \times 5$ image $V$ with a $3 \times 3$ convolutional operator $U$ that employs one filter (i.e., $q = 1$).

\begin{figure}[!htb]
    \centering
    \includegraphics[width=0.5\textwidth]{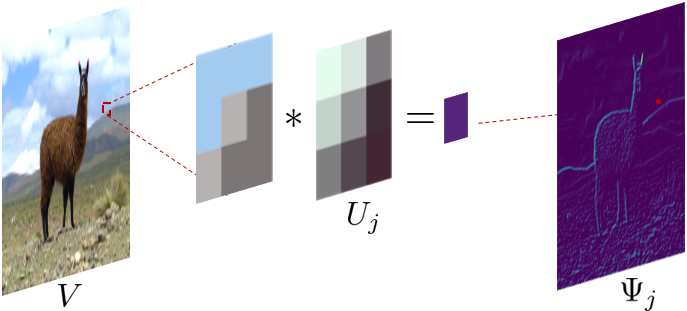}
    \caption{Illustration of how the filter $U_j$ is convolved with $V$ over its neighborhoods to produce feature map matrix $\Psi_j$ which shows how the pattern in $U_j$ is manifested in $V$.}
    \label{fig:flare_map}
\end{figure}

Following Figure \ref{fig:convolution}, a common and intuitive interpretation of applying convolutional operators is in terms of pattern recognition. The $q$ filters that comprise $U$ each embed a particular pattern that when convolved with a particular neighborhood of $V$ (as depicted in Figure \ref{fig:convolution}) gives a score of how well the pattern is matched (where larger values denote greater matching). Hence, the feature map $\Psi$ can be interpreted as a score sheet of how certain patterns (as encoded by the filters in $U$) are manifested in $V$. Here, we let $U_j$ denote the $j^\text{th}$ convolutional filter and $\Psi_j$ denotes the feature map slice (a matrix) that corresponds to $U_j$. Figure \ref{fig:flare_map} demonstrates how convolving an image with a filter that encodes an edge pattern gives a feature map that records the highest scores in the regions of the image that have similar edges.

The feature map output of a convolutional layer is typically mapped element-wise through an activation function $\alpha : \mathbb{R} \mapsto \mathbb{R}$ to yield the activation $A \in \mathbb{R}^{n_\psi \times n_\psi \times q}$:
\begin{equation}
    A_{x_1,x_2,j} = \alpha(\Psi_{x_1,x_2,j}).
    \label{eq:activation}
\end{equation}
The functional form of this operation employs $f_a : \mathbb{R}^{n_\psi \times n_\psi \times q} \mapsto \mathbb{R}^{n_\psi \times n_\psi \times q}$ such that $A = f_a(\Psi)$. This operation helps the CNN sensor to encode nonlinear behavior. Common choices of $\alpha(\cdot)$ include:
\begin{equation*}
    \begin{gathered}
        \alpha_\text{sig}(z) = \frac{1}{1 + e^{-z}} \\
        \alpha_\text{tanh}(z) = \text{tanh}(z) \\ 
        \alpha_\text{ReLU}(z) = \text{max}(0, z).
    \end{gathered}
\end{equation*}
Here the Rectified Linear Unit (ReLU) function $\alpha_\text{ReLU}(\cdot)$ has achieved heightened popularity since it generally exhibits greater sensitivity to changes in the input \cite{nair2010rectified}.

\begin{figure}[!htb]
    \centering
    \includegraphics[width=0.4\textwidth]{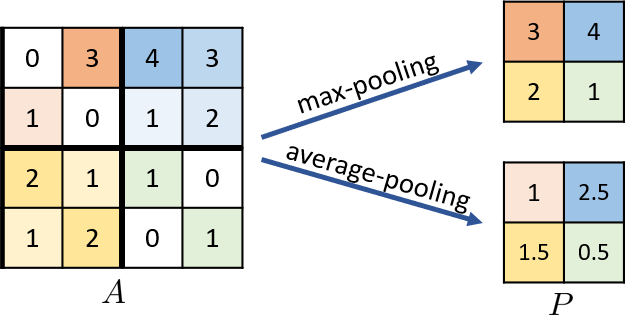}
    \caption{An illustration of max-pooling and average-pooling.}
    \label{fig:pooling}
\end{figure}

Another key component of CNNs are pooling layers. Pooling operations are dimension-reduction mappings $f_p : \mathbb{R}^{n_\psi \times n_\psi \times q} \mapsto \mathbb{R}^{n_P \times n_P \times q}$ that seek to summarize/reduce the activation $A$ by collapsing certain sub-regions of dimensions $n_p \times n_p$ to scalar values (referred to as pooling). Here we have that $n_P = n_\psi / n_p$. The functional form this operation is expressed $P = f_p(A)$ with the pooled output signal $P \in \mathbb{R}^{n_P \times n_P \times q}$. Pooling helps make the learned representation more invariant to small length-scale perturbations \cite{nagi2011max}. Common choices include max-pooling and average-pooling where the maximum or average of a sub-region is used to scalarize, respectively. Figure \ref{fig:pooling} illustrates these typical operations.

\begin{figure}[!htb]
     \centering
     \begin{subfigure}[b]{0.3\textwidth}
         \centering
         \includegraphics[width=\textwidth]{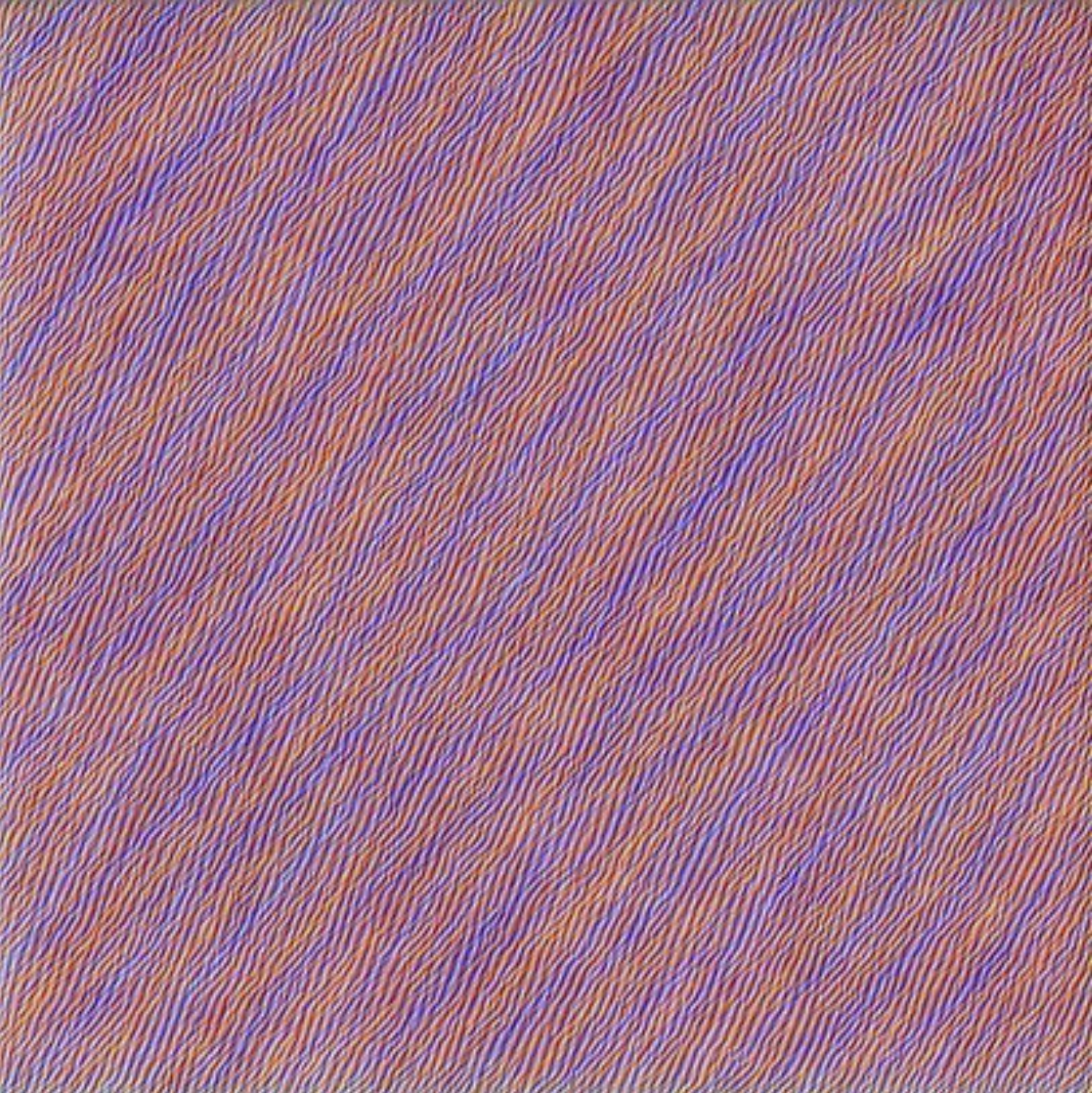}
         \caption{Block 1, Filter 16}
     \end{subfigure}
     \quad
     \begin{subfigure}[b]{0.3\textwidth}
         \centering
         \includegraphics[width=\textwidth]{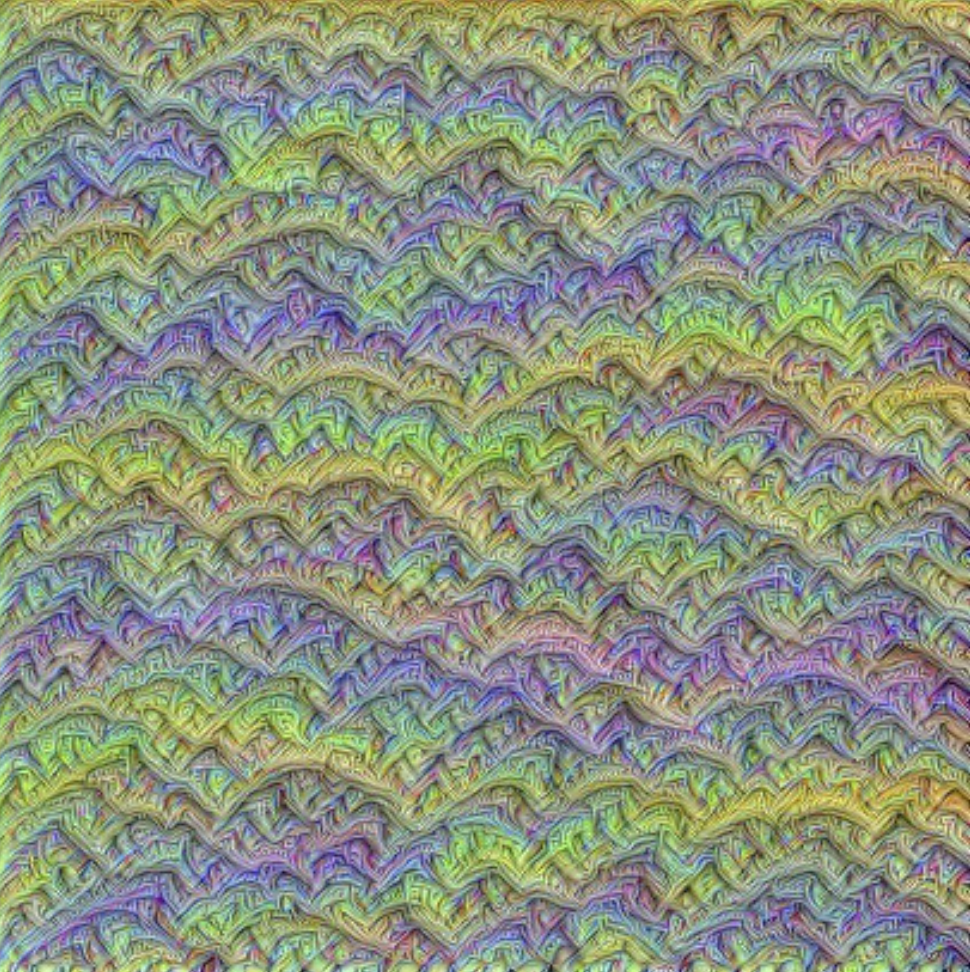}
         \caption{Block 3, Filter 181}
     \end{subfigure}
     \quad
     \begin{subfigure}[b]{0.3\textwidth}
         \centering
         \includegraphics[width=\textwidth]{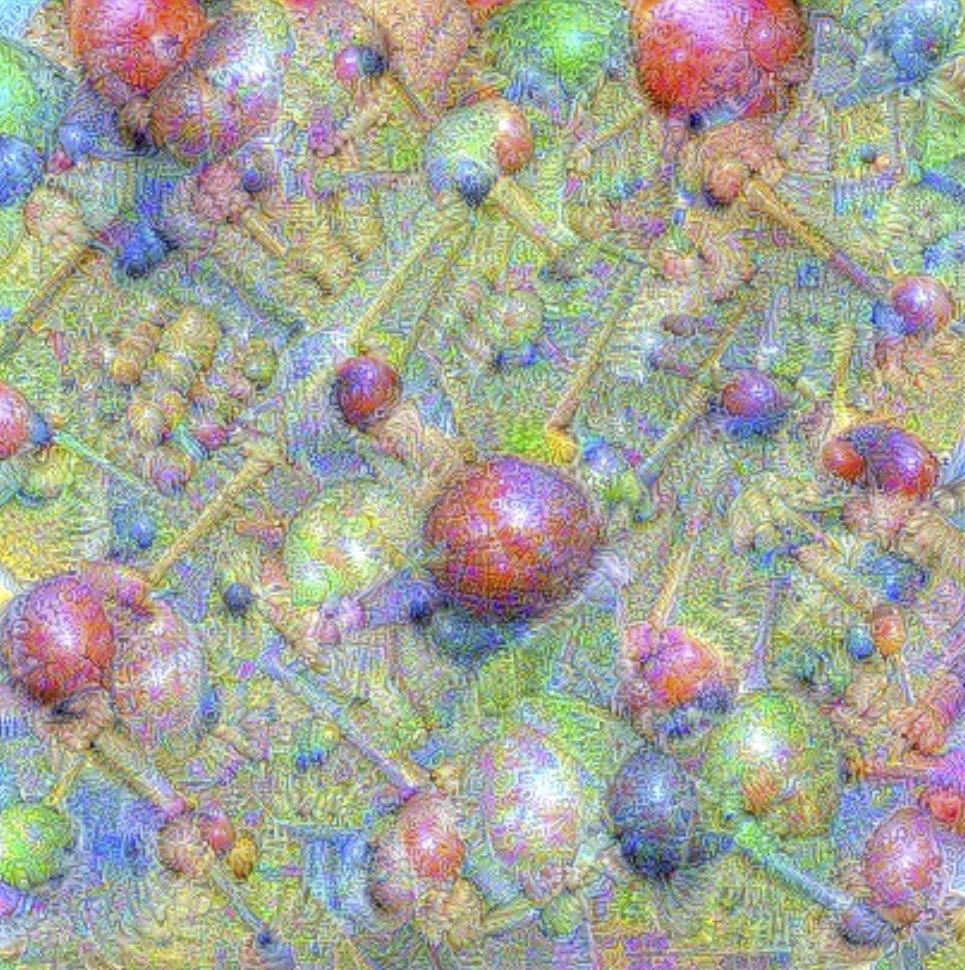}
         \caption{Block 5, Filter 261}
     \end{subfigure}
    \caption{The images that maximize the mean output of the particular filters in the VGG16 CNN model. The filters from deeper blocks detect more refined patterns/objects.}
    \label{fig:max_features}
\end{figure}

Convolutional blocks in CNNs are a combination of convolution, activation, and pooling a given input. In other words, a convolution block employs the mapping $f_{\text{cb}} : \mathbb{R}^{n_v \times n_v \times p} \mapsto \mathbb{R}^{n_P \times n_P \times q}$ such that we have:
\begin{equation}
    P = f_{\text{cb}}(V; U) = f_p(f_a(f_c(V; U))).
    \label{eq:convolution_block}
\end{equation}
Note that blocks can employ more complex operation nesting (e.g., multiple convolutional layers), but we consider blocks that follow \eqref{eq:convolution_block} for simplicity in presentation. Moreover, these can take pooled outputs $P$ as input and thus facilitate the use of multiple convolutional blocks in succession via recursively calling $f_{\text{cb}}$. We will identify the block that a particular function or tensor belongs to in a CNN via the superscript index ${(\ell)}$ where $\ell \in \{1, n\}$. Hence, for a CNN that employs a couple of convolutional blocks, we have:
\begin{equation*}
    \begin{gathered}
        P^{(1)} = f_{\text{cb}}^{(1)}(V; U^{(1)}) \\ 
        P^{(2)} = f_{\text{cb}}^{(2)}(P^{(1)}; U^{(2)}).
    \end{gathered}
\end{equation*}
These blocks act as feature extractors via the $q$ convolutional filters they employ. Moreover, this feature extraction becomes more specialized for blocks that are located deeper in the CNN. In the context of computer vision, this means that the first blocks extract simple features (e.g., edges or colors) and the deeper blocks can extract more complex patterns (e.g., particular and abstract shapes). We can think of performing multiple convolution block operations as a multi-stage distillation process where the successive feature spaces capture patterns of increased length-scale and complexity (i.e., distill the image data into an increasingly purified feature space) \cite{chollet2017deep}. Figure \ref{fig:max_features} demonstrates this principle by showing the images that maximally activate a particular filter in the first, third, and fifth convolutional blocks of the popular VGG16 CNN model. We will show how these types of convolutional blocks are leveraged in the SAFE-OCC framework to extract effective feature spaces for novelty detection.

\begin{figure}[!htb]
    \centering
    \includegraphics[width=0.9\textwidth]{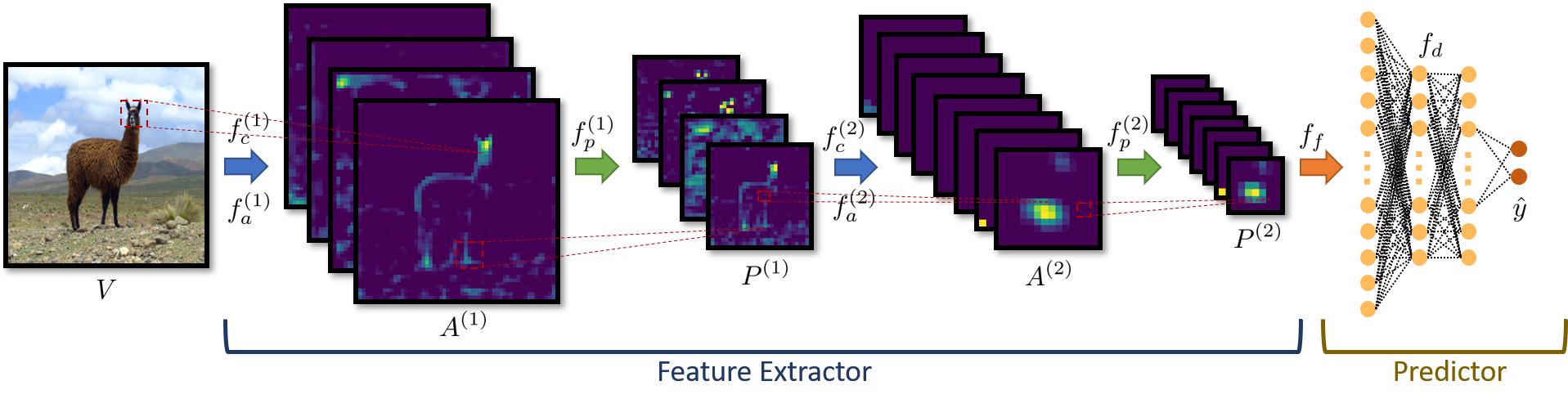}
    \caption{Schematic of a typical CNN sensor.}
    \label{fig:cnn_example}
\end{figure}

The output $P$ of the last convolutional block is typically flattened in the mapping $f_f : \mathbb{R}^{n_P \times n_P \times q} \mapsto \mathbb{R}^{n_P \cdot n_P \cdot q}$ to yield the feature vector $v \in \mathbb{R}^{n_P \cdot n_P \cdot q}$:
\begin{equation}
    v = f_f(P).
\end{equation}
The feature vector is then fed into a dense neural network model $f_d : \mathbb{R}^{n_P \cdot n_P \cdot q} \mapsto \mathbb{R}^{n_y}$ which predicts the desired state space vector $\hat{y}$ (often a regression problem in the context of process control). Figure \ref{fig:cnn_example} illustrates a typical CNN model that implements the components described above. It employs two convolutional blocks and can be described in the functional form:
\begin{equation}
    \hat{y} = f_{\text{cnn}}(V) = f_d\left(f_f\left(f_{\text{cb}}^{(2)}\left(f_{\text{cb}}^{(1)}\left(V; U^{(1)}\right); U^{(2)}\right)\right)\right).
\end{equation}
This again emphasizes that the convolutional blocks act as feature extractors and the dense layers act as the predictor whose feature space is the flattened output $v$ of the final convolutional block.

Training procedures seek optimal model parameters (i.e., convolution operators, dense network weights) that minimize the error incurred by the state predictions made relative to the training data set. Here, we consider a training set $\{(V^{(k)}, y^{(k)}) : k \in \mathcal{K}\}$ that employs $|\mathcal{K}|$ image-state pairs. The prediction error minimized in the training procedure is called the loss function $L : \mathbb{R}^{n_y} \mapsto \mathbb{R}$. For example, regression models typically use a sum-of-squared-error (SSE) loss function:
\begin{equation}
    L(\hat{y}) = ||\hat{y} - y||_2^2.
\end{equation}
Thus, by grouping all the CNN model parameters into $\theta \in \mathbb{R}^{n_\theta}$ we can express model training as a standard optimization problem:
\begin{equation}
    \begin{aligned}
        &\min_\theta &&  \sum_{k \in \mathcal{K}} L(\hat{y}^{(k)}) \\
        &\text{s.t.} && \hat{y}^{(k)} = f_{\text{cnn}}(V^{(k)}; \theta), && k \in \mathcal{K}.
    \end{aligned}
\end{equation}
Note that this can readily be expressed as an unconstrained optimization problem by inserting the constraint equations directly into the objective. Stochastic Gradient Descent (SGD) is typically used to solve this problem due to the large amount of training data, the high number of model parameters, and the model complexity. Moreover, forward and backward propagation techniques are used to evaluate the objective and derivative values required by each iteration of the SGD algorithm.
\\

\begin{figure}[!htb]
     \centering
     \begin{subfigure}[b]{0.2\textwidth}
         \centering
         \includegraphics[width=\textwidth]{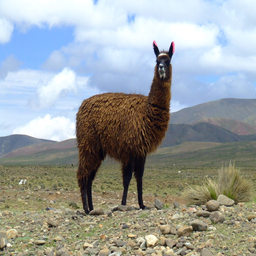}
         \caption{Original}
     \end{subfigure}
     \quad
     \begin{subfigure}[b]{0.2\textwidth}
         \centering
         \includegraphics[width=\textwidth]{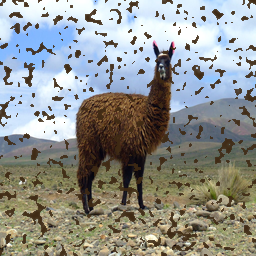}
         \caption{Splattered}
     \end{subfigure}
     \quad
     \begin{subfigure}[b]{0.2\textwidth}
         \centering
         \includegraphics[width=\textwidth]{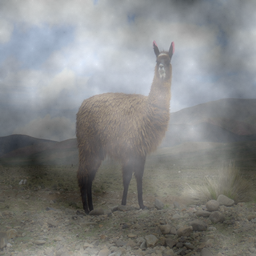}
         \caption{Fogged}
     \end{subfigure}
     \quad
     \begin{subfigure}[b]{0.2\textwidth}
         \centering
         \includegraphics[width=\textwidth]{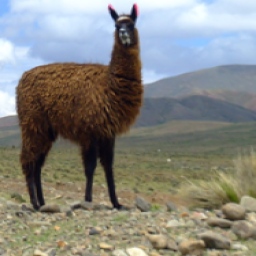}
         \caption{Shifted}
     \end{subfigure}
    \caption{Examples of image augmentation.}
    \label{fig:image_aug}
\end{figure}

{\em Image augmentation} is often used to expand the size of the training image set in an effort to decrease the likelihood of the CNN sensor encountering novel images. Image augmentation generally denotes perturbing the training images such that the CNN sensor can be robust to those types of visual disturbance. Common perturbations include rotation, translation, cropping, blurring, brightness changing, splattering, and more. There are many software tools available to implement these transformations which include \texttt{TensorFlow} and \texttt{ImgAug} \cite{geron2019hands, imgaug}. Figure \ref{fig:image_aug} exemplifies how a training image is augmented via a variety of perturbation (disturbance) types. This methodology helps mitigate the risk of CNN sensors encountering novel images, but it is not typically possible to account for all the disturbance a process might encounter.

We will show how the SAFE-OCC framework can be readily incorporated with existing CNN sensors to viably detect novel images relative to the feature space of the CNN. Moreover, the addition of the SAFE-OCC framework to a workflow typically incurs only a minor increase in computational cost.

\subsection{PCA-Based Dimension Reduction} \label{sec:pca}

\begin{figure}[!htb]
    \centering
    \includegraphics[width=0.4\textwidth]{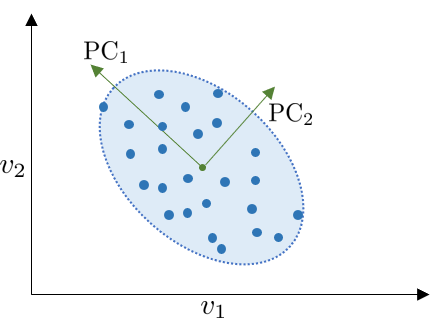}
    \caption{An illustration of two-feature data points and their corresponding principal components.}
    \label{fig:1dpca}
\end{figure}

Principal Component Analysis (PCA) is a popular dimensionality reduction technique for 1D vector data. It achieves this reduction by projecting an input vector $v \in \mathbb{R}^{n_v}$ onto the column-space $W \in \mathbb{R}^{n_v \times d}$:
\begin{equation}
    f_\text{pca}(v)^T = v^TW
\end{equation}
where we have that $d \leq n_v$ and $f_\text{pca} : \mathbb{R}^{n_v} \mapsto \mathbb{R}^{d}$ is the PCA mapping function. The columns of the linear projection matrix $W$ are derived from the eigenvectors that correspond to the $d$ largest eigenvalues of the empirical covariance matrix $\Sigma \in \mathbb{R}^{n_v \times n_v}$. Here $\Sigma$ can be computed via the outer product of the mean centered data:
\begin{equation}
    \Sigma = \frac{1}{|\mathcal{K}|} \sum_{k \in \mathcal{K}} (v^{(k)} - \bar{v}) (v^{(k)} - \bar{v})^T
\end{equation}
where $\bar{v} := |\mathcal{K}|^{-1} \sum_{k \in \mathcal{K}} v^{(k)}$ is the empirical average over the training set $\{v^{(k)} : k \in \mathcal{K}\}$. Hence, PCA projects 1D data into a reduced space whose bases are the so-called principal components. The orthogonal principal components capture the directions of maximum variance in the training data since they correspond to the principal axes of the ellipsoid with weight matrix $\Sigma$ as illustrated in Figure \ref{fig:1dpca}. Moreover, the amount of variance captured by a particular principal component is defined by its eigenvalue (i.e., the principal components with the largest eigenvalues describe the directions of greatest variance). 

High-dimensional data (e.g., images and feature maps) are not well-suited for PCA since flattening them into a 1D data representation will result in having  numerous features relative to the number of training samples; moreover, flattening can eliminate information on spatial context and correlation. To address this, Zhang and Zhou proposed a Two-Directional Two-Dimensional Principal Component Analysis (2D\textsuperscript{2}PCA) in \cite{zhang20052d} to generalize PCA for 2D matrix data and to address these concerns. This approach builds upon the 2DPCA approach proposed in \cite{yang2004two}, but 2D\textsuperscript{2}PCA differs in its ability to reduce matrix data over both rows and columns where the 2DPCA is only able to reduce over the rows of a matrix. Here, we reduce a 2D matrix $V \in \mathbb{R}^{n_v \times n_v}$ (square for simplicity in presentation) via two projection matrices $W \in \mathbb{R}^{n_v \times n_d}$ and $Q \in \mathbb{R}^{n_v \times n_r}$:
\begin{equation} 
    f_\text{2d2pca}(V) = Q^TVW
    \label{eq:2d2pca}
\end{equation}
where we have that $d, r \leq n_v$ and $f_\text{2d2pca} : \mathbb{R}^{n_v \times n_v} \mapsto \mathbb{R}^{r \times d}$ is the 2D\textsuperscript{2}PCA mapping function. The projection matrices $W$ and $Q$ reduce the matrix $V$ over its columns and rows, respectively. These are learned from covariance matrices $\Sigma_W \in \mathbb{R}^{n_v \times n_v}$ and  $\Sigma_Q \in \mathbb{R}^{n_v \times n_v}$ which are defined:
\begin{equation}
    \begin{aligned}
        \Sigma_W &= \frac{1}{|\mathcal{K}|} \sum_{k \in \mathcal{K}} (V^{(k)} - \bar{V})^T (V^{(k)} - \bar{V}) \\ 
        \Sigma_Q &= \frac{1}{|\mathcal{K}|} \sum_{k \in \mathcal{K}} (V^{(k)} - \bar{V}) (V^{(k)} - \bar{V})^T
    \end{aligned}
\end{equation}
where $\bar{V} := |\mathcal{K}|^{-1} \sum_{k \in \mathcal{K}} V^{(k)}$ is the empirical average over the training set $\{V^{(k)} : k \in \mathcal{K}\}$. We then derive the columns of $W$ and $Q$ using the eigenvalues and eigenvectors of their respective covariance matrices in like manner to PCA. In comparison to using traditional PCA with flattened matrix data, Zhang and Zhou find that 2D\textsuperscript{2}PCA is able to derive a feature a space that leads to better classification accuracy and incurs a lower computational cost in the context of facial recognition. We will demonstrate that this approach can effectively scalarize CNN feature maps (matrices) within the SAFE-OCC framework.

\subsection{One-Class Classification} \label{sec:occ} 

One-class classification denotes a group of novelty detection approaches that identify a single class (group) of normal data (based on training data with no or few novel instances) that then is used to discriminate whether new data is normal or novel. OCC is widely applied and thus many such approaches have been developed. Examples include One-Class Support Vector Machines (OC-SVMs), One-Class Classifier Ensembles, Neural Network Models, Decision Trees, Bayesian Classifiers, and more. Moreover, many of these can be readily implemented with software tools such as \texttt{SciKit-Learn} \cite{pedregosa2011scikit}. Reviewing all of these approaches is beyond the scope of this work, but we will generally describe the functionality of OCC in our context using OC-SVMs as a concrete example. Here, the input data point $v \in \mathbb{R}^{n_v}$ contains $n_v$ features. The OCC method is represented via the mapping $f_{\text{occ}} : \mathbb{R}^{n_v} \mapsto \mathbb{R}$ which gives the predicted output $\hat{h} \in \mathbb{R}$. Typically, the output $\hat{h}$ is thresholded by $\rho \in \mathbb{R}$ to determine if $v$ is normal or novel. Moreover, training involves selecting the optimal model parameters to identify normal samples. 

\begin{figure}[!htb]
    \centering
    \includegraphics[width=0.7\textwidth]{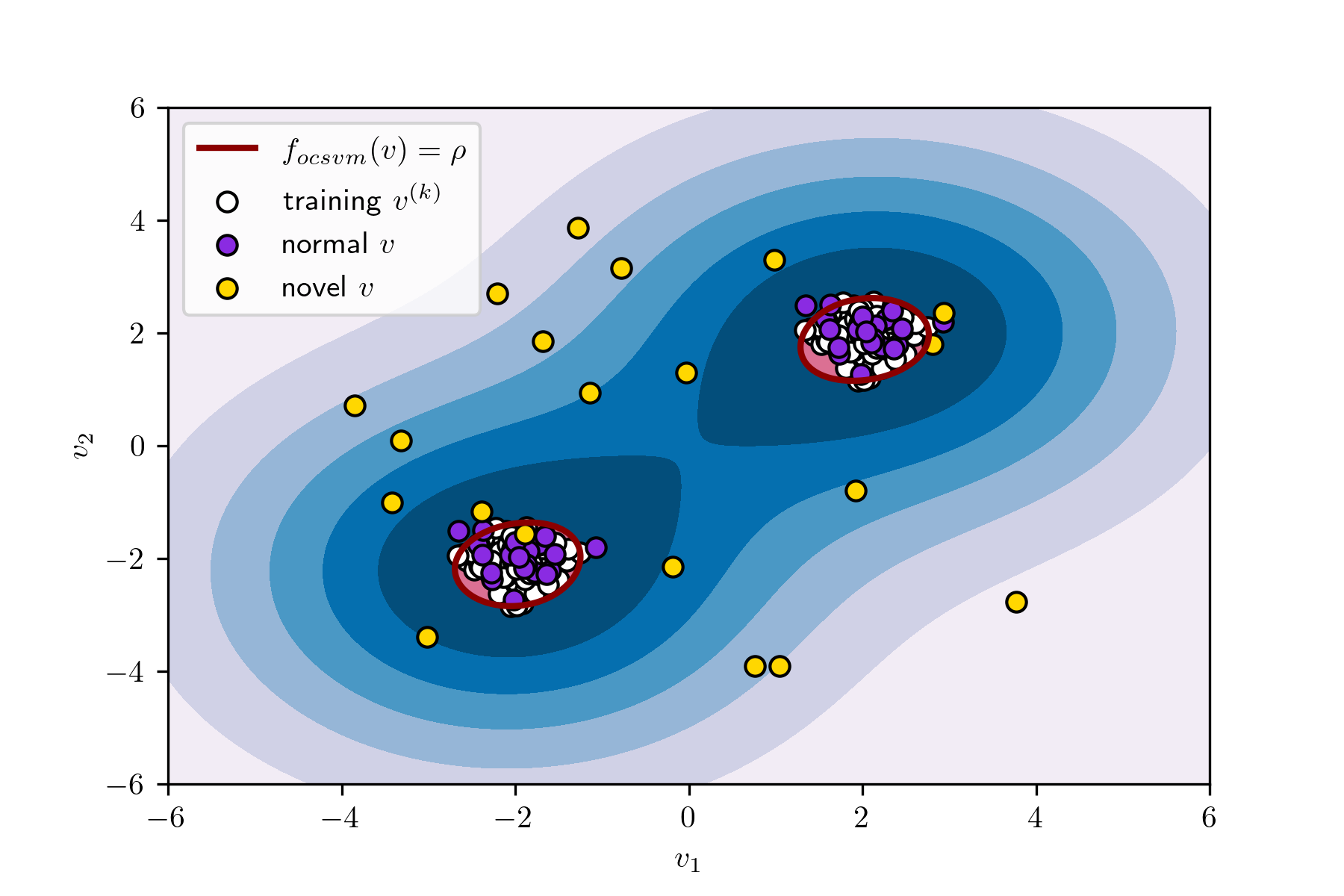}
    \caption{Learned threshold boundary $f_{\text{ocsvm}}(v) = \rho$ for a OC-SVM in a simple feature space. We observe that the novelty of a point is determined by its Euclidean distance to the learned boundary.}
    \label{fig:ocsvm_example}
\end{figure}

OC-SVMs typically learn a boundary around the training data to a certain threshold whose shape is determined by a kernel function $\kappa(v^{(i)}, v^{(j)})$ which quantifies the similarity between samples. One common choice is the Gaussian kernel:
\begin{equation}
    \kappa(v^{(i)}, v^{(j)}) = \exp\left(-\gamma||v^{(i)} - v^{(j)}||^2_2\right)
\end{equation}
where $\sigma \in \mathbb{R}_+$. Thus, the mapping function $f_{\text{ocsvm}}$ is:
\begin{equation}
    \hat{h} = f_{\text{ocsvm}}(v) = \sum_{k \in \mathcal{K}_{sv}} \alpha_k \kappa(v, v^{(k)})
\end{equation}
where $\mathcal{K}_{sv}$ is the set of indices of support vectors (taken from the training data) and $\alpha_k \in \mathbb{R}$ are Lagrangian coefficients that satisfy $\sum_{k \in \mathcal{K}_{sv}} \alpha_k = 1$ \cite{bounsiar2014kernels}. Here, a particular sample is novel if $\hat{h} < \rho$, otherwise it is classified as normal. Training involves choosing support vectors $v^{(k)}$ (with the corresponding values of $\alpha_k$) and selecting the threshold parameter $\rho$ to derive a boundary that well encloses the training data. Figure \ref{fig:ocsvm_example} illustrates a Gaussian kernel OC-SVM for a two-feature system. We will show how OCC approaches can be readily incorporated into the SAFE-OCC framework to effectively assess the performance of a CNN sensor. 

\subsection{Process Control} \label{sec:control_back}

Process control is a broad discipline for automation systems that manipulate the inputs of a process to achieve/maintain desired operation states. In this work, we consider a process with control input variables $z \in \mathbb{R}^{n_z}$ and state output variables $y \in \mathbb{R}^{n_y}$. Typically, measurement sensors (i.e., thermocouples, pressure transducers, flow meters, speedometers) are used to get a measured estimate $\hat{y} \in \mathbb{R}^{n_y}$ of the state variables $y$. These measurements are compared against the setpoint $y_\text{sp} \in \mathbb{R}^{n_y}$ (the desired operation state) to yield the state error $y_e \in \mathbb{R}^{n_y}$ (i.e., the setpoint tracking error). This setpoint tracking error is the input to a controller which outputs the control variables $z$. 
\\

\begin{figure}[!htb]
    \centering
    \includegraphics[width=0.7\textwidth]{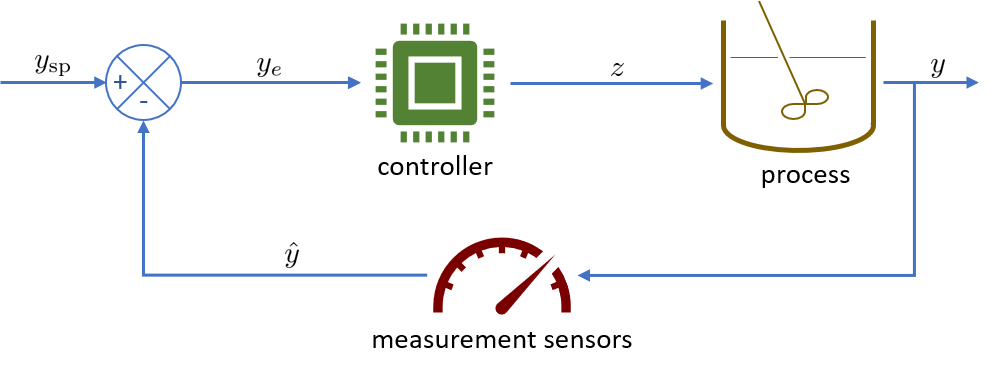}
    \caption{A general feedback control loop schematic.}
    \label{fig:feedback_control}
\end{figure}

Diverse control system configurations are used in practice in accordance with the characteristics of the process being automated. Common paradigms include feedback, feed-forward, and cascade control. For simplicity in presentation, we consider feedback control loops as illustrated in Figure \ref{fig:feedback_control}. This system features the procedure where control variables $z$ are given to a process whose state variable measurements $\hat{y}$ are compared against the setpoint $y_\text{sp}$ to determine the error $y_e$ which is given to the controller that determines the next control variables $z$. The procedure is conducted continuously to automatically control the process to follow the desired operational setpoint. Typical control paradigms are proportional-integral-derivative (PID) control and model predictive control (MPC). For example, a PID controller has the control law:
\begin{equation}
    z(t) = K_p y_e(t) + K_i \int y_e(t) dt + K_p \frac{d y_e(t)}{dt}
    \label{eq:pid}
\end{equation}
where $K_p,K_i \in \mathbb{R}$ are constant parameters and the variables $z$ and $y_e$ are shown as time-valued functions to index their values over time $t$.

Control systems generally consider configurations of greater complexity that incorporate a variety of other functions such as disturbance estimation, Kalman filtering, and more (which are beyond the scope of this work). In this work, our focus will be targeted on the implications of using CNN sensors as nontraditional measurement/sensing devices that can convert high-dimensional image data into state signals that can be readily incorporated into existing process control architectures. Moreover, we will see how this computer vision aided process control paradigm motivates the development of our SAFE-OCC novelty detection framework.

\section{CNN Sensor Aided Process Control} \label{sec:compvis_control}

In this section, we build upon the discussion in Section \ref{sec:control_back} to discuss the motivation and implications for incorporating CNN sensors into process control systems. We consider processes that have state variables $y_\text{trad} \in \mathbb{R}^{n_{yt}}$ which are measured via conventional sensors and state variables $y_\text{vis} \in \mathbb{R}^{n_{yv}}$ which can be measured from visual observation (e.g., a camera feed). These can be concatenated to yield the full vector of state variables $y$ (i.e., $n_y = n_{yt} + n_{yv}$).
\\

\begin{figure}[!htb]
    \centering
    \includegraphics[width=0.8\textwidth]{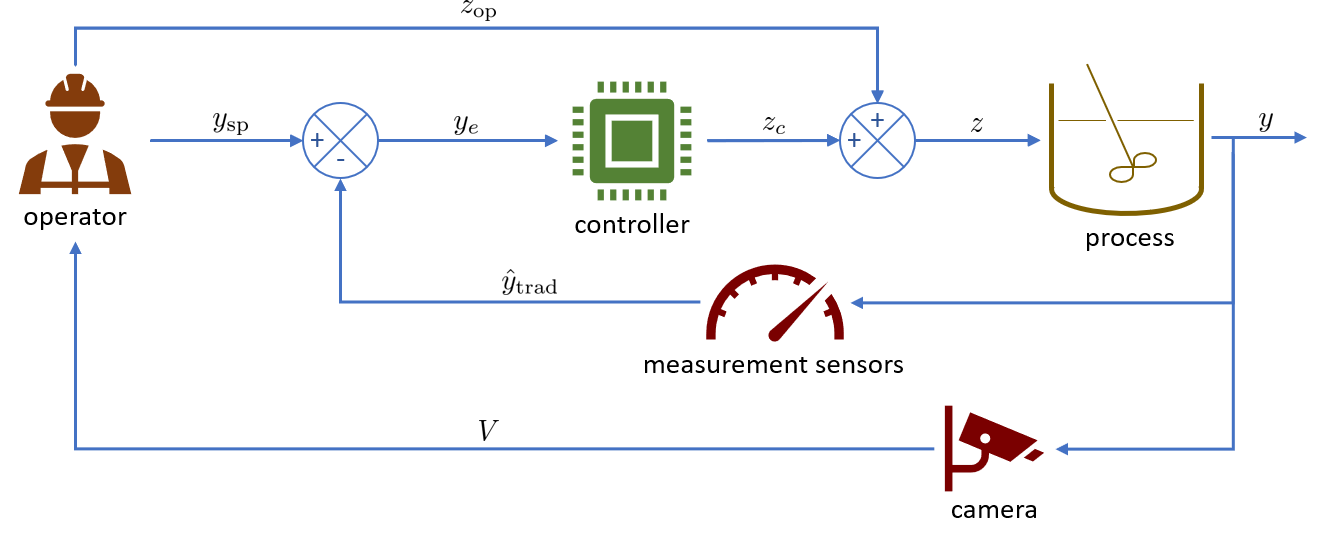}
    \caption{Feedback control loop where photogrammetry data is used by a human operator to make control/operation decisions while the process is simultaneously controlled automatically using traditional measurement signals. Here, the operator becomes incorporated into the control loop.}
    \label{fig:operator_control}
\end{figure}

Figure \ref{fig:operator_control} shows a traditional feedback control system where an automatic control loop operates using traditional measurement sensors. Meanwhile, an operator will monitor the visual data $V$ which he/she will implicitly use to predict the values of $y_\text{vis}$ and make adjustments to the setpoint $y_\text{sp}$ and/or intervene with manual control action $z_\text{op} \in \mathbb{R}^{n_z}$ (which is combined with the automatic control variables $z_c \in \mathbb{R}^{n_z}$) as needed. Here the operator effectively becomes part of the control loop, making this control system an intermediate between closed-loop and open-loop control. This paradigm can be quite burdensome for an operator that typically has to monitor and control many systems simultaneously. An everyday example would be a driver that uses visual cues to adjust the setpoint of a vehicle cruise control and/or temporarily press the throttle peddle to an extra extent. Industrial examples could include flare visible emission control and controlling the flow of polymer material through an extruder. 

\begin{figure}[!htb]
    \centering
    \includegraphics[width=0.8\textwidth]{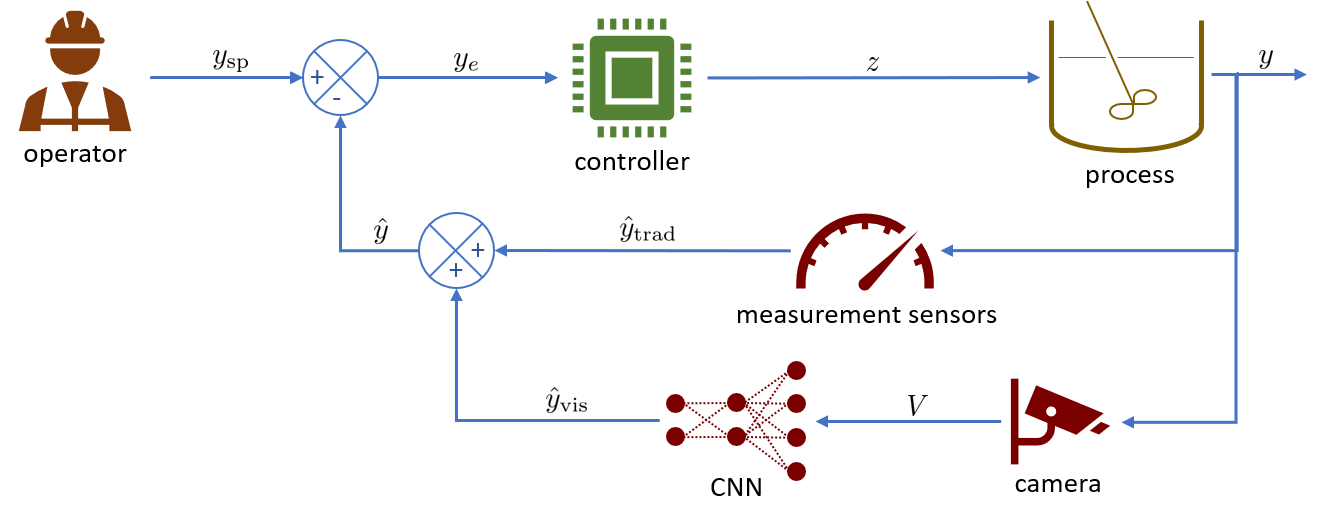}
    \caption{Feedback control loop that incorporates a CNN to convert photogrammetry data into a measurement signal that is amendable for automatic control; thus, negating then need for an operator to manually interpret/control it.}
    \label{fig:vision_control}
\end{figure}

We consider leveraging a CNN sensor to autonomously map image data $V$ to visual state variables $y_\text{vis}$ such that we achieve closed-loop control. With this, we remove the operator from the control loop in the sense that he/she will no longer need to actively interpret and act upon visual process data. Such a control system is depicted in Figure \ref{fig:vision_control}. Here, the camera and the CNN work together to form a computer vision sensor that is able to measure states $y_\text{vis}$ that otherwise would not be available using traditional process measurement devices. Hence, following this new paradigm we obtain a fully automatic control system that can exhibit improved setpoint tracking performance which promotes increased process safety, profitability, and consistency. 

The paradigm shift from an operator-centric control system of Figure \ref{fig:operator_control} to the CNN-aided system of Figure \ref{fig:vision_control} introduces a significant vulnerability: poor prediction accuracy of $y_\text{vis}$ when $V$ is novel relative to the training data used to prepare the CNN (i.e., the CNN sensor makes a highly inaccurate prediction because it is extrapolating). Injecting erroneous measurement data into a closed-loop control architecture can have severe consequences to profitability and safety. Image augmentation, as described in Section \ref{sec:cnn}, can be used to help alleviate this problem by seeking to account for a variety of visual disturbances a process might encounter (thus increasing the span of the training image set). Adversarial training can also be used to perturb training data and enhance robustness. These approaches can help make extrapolation events more rare, but it is usually not possible to account for every possible visual disturbance a process might be subjected to. Thus, we require an approach to automatically recognize when the visual data $V$ is novel relative to the CNN sensor being used. Such an approach can be incorporated into a monitoring/safety system that will mitigate the risk of unknowingly injecting inaccurate CNN sensor measurements into a control system. This need is what motivates the creation of the SAFE-OCC novelty detection framework which we propose and detail in Section \ref{sec:framework}.

\section{SAFE-OCC Novelty Detection Framework} \label{sec:framework}

In this section, we detail the SAFE-OCC novelty detection framework that leverages the native feature space of a CNN sensor to achieve novelty detection that is complimentary. This framework is principally comprised of three steps: feature extraction via the feature maps of a CNN sensor, feature refinement, and novelty detection via OCC. We will see that this framework creates OCC novelty detectors whose feature space closely relates to  that of the targeted CNN sensor, in contrast to conventional approaches that derive a feature space independently (potentially making them less effective in identifying novel data relative to the CNN sensor). It is this distinction that makes this framework a natural fit for CNN sensors in process control applications.

\subsection{Feature Extraction} \label{sec:feature_extract}

We recall that CNNs use convolutional blocks $f_\text{cb}$ to extract features from high-dimensional image data $V$, which is not readily amendable for the dense network predictor $f_d$, to ultimately derive a 1D feature space (i.e., comprised of vector data points $v$). Similarly, OCC methods typically require a 1D input $v$ and thus are not directly compatible with high-dimensional image data. Since we are interested in assessing the novelty of image data relative to the predictor $f_d$, we would like the OCC novelty detector to use a similar input feature space (not independently derive one from the image data). Naively, we could directly use the feature vector output $v$ of the flattening mapping $f_f$, but this will typically have a large number (e.g., thousands) of features which can be prohibitively large for training effective OCC novelty detectors. 

\begin{figure}[!htb]
    \centering
    \includegraphics[width=0.8\textwidth]{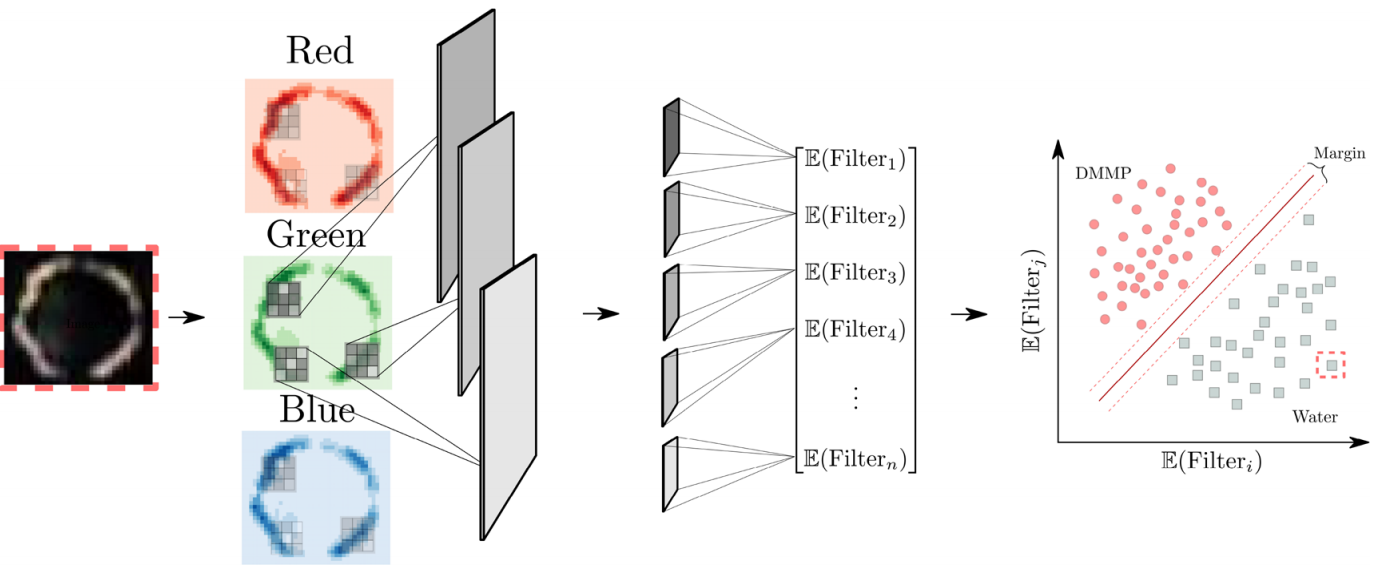}
    \caption{Illustration of feature extraction framework to distinguish the exposure of liquid crystal sensors to DMMP vs. water  \cite{smith2020convolutional}. Here, $\mathbb{E}(\cdot)$ denotes the spatial average used to scalarize each feature map matrix $P_j$.}
    \label{fig:lc_sensor}
\end{figure}

In the proposed approach, we seek a reduced feature space that is derived intelligently from the feature map output $P^{(\ell)}$ of a particular convolutional block $f_\text{cb}^{(\ell)}$. Strictly speaking, $P$ denotes the output of feature maps that have been activated and pooled, but with some abuse of nomenclature we refer to these as feature maps for convenience in presentation. In \cite{andrews2016transfer}, Andrews and colleagues use the output $P$ of the terminal layer of the pre-trained VGG-F and VGG-M CNN models to derive a feature space for novelty detection in an approach they call Transfer Representation-Learning Anomaly Detection. They apply this approach to multiple datasets (none of which were part of the CNN training data) and generally achieve viable classification accuracy. Similarly, in \cite{smith2020convolutional}, Smith and colleagues derive a reduced feature space from the output of the first convolutional block of the VGG16 CNN model. They observe that a tensor feature map output $P \in \mathbb{R}^{n_P \times n_P \times q}$ can be broken up into a set of $q$ feature map matrices $P_j \in \mathbb{R}^{n_P \times n_P}$ that each corresponds to a convolutional filter $U_j$ (as illustrated in Figures \ref{fig:flare_map} and \ref{fig:cnn_example}); thus, they derive their feature space by scalarizing each feature map matrix $P_j$ via its average to obtain a vector of scalarized features. With this methodology, they produce an effective feature space for liquid crystal sensor images that yields high accuracy classification with the SVM used in their study. This approach is summarized in Figure \ref{fig:lc_sensor}.
\\

Taking inspiration from these approaches, we propose the scalarizing operator $f_s : \mathbb{R}^{n_P \times n_P \times q} \mapsto \mathbb{R}^q$ that derives a feature space with data points $v \in \mathbb{R}^q$. This applies a scalarizing function $g : \mathbb{R}^{n_P \times n_P} \mapsto \mathbb{R}$ to each of the $q$ feature map matrices $P_j$. Hence, the reduced representation $v$ comprises $q$ elements that each summarize the feature map output of a convolutional filter that assesses the prominence of a certain pattern in the input image $V$. Candidate choices of $g(\cdot)$ include:
\begin{equation*}
    \begin{gathered}
        g_\text{max}(P_j) = \max(P_j) \\ 
        g_\mathbb{E}(P_j) = \mathbb{E}(P_j)\\ 
        g_\text{2d2pca}(P_j) = f_\text{2d2pca}(P_j)
    \end{gathered}
\end{equation*}
where $\max(\cdot)$ returns the largest element of an input matrix, $\mathbb{E}(\cdot)$ returns the average value of a matrix, and $f_\text{2d2pca}(\cdot)$ reduces a matrix to a scalar value via 2D\textsuperscript{2}PCA. The scalarization $f_s$ can be interpreted as a special case of the pooling function $f_p$ where a single pooling operation is done over the entirety of each feature map matrix $P_j$ that comprises $P$. With this interpretation, the maximization $\max(\cdot)$ captures the most activated presence of a feature (as assessed by a convolutional filter) and the average $\mathbb{E}(\cdot)$ summarizes the presence of a feature \cite{nagi2011max}. The use of 2D\textsuperscript{2}PCA on feature map matrices has not been explored in the literature to the best of our knowledge, however, the results presented below in Section \ref{sec:case1} suggests that it can be quite effective at scalarizing/summarizing feature maps. This could be attributed to the ability of 2D\textsuperscript{2}PCA to derive reduced representations that incorporate spatial relationships present in the feature maps. Note that the 2D\textsuperscript{2}PCA model should be trained using feature map matrices that derive from the same training image data used by the CNN sensor.

\begin{figure}[!htb]
    \centering
    \includegraphics[width=0.9\textwidth]{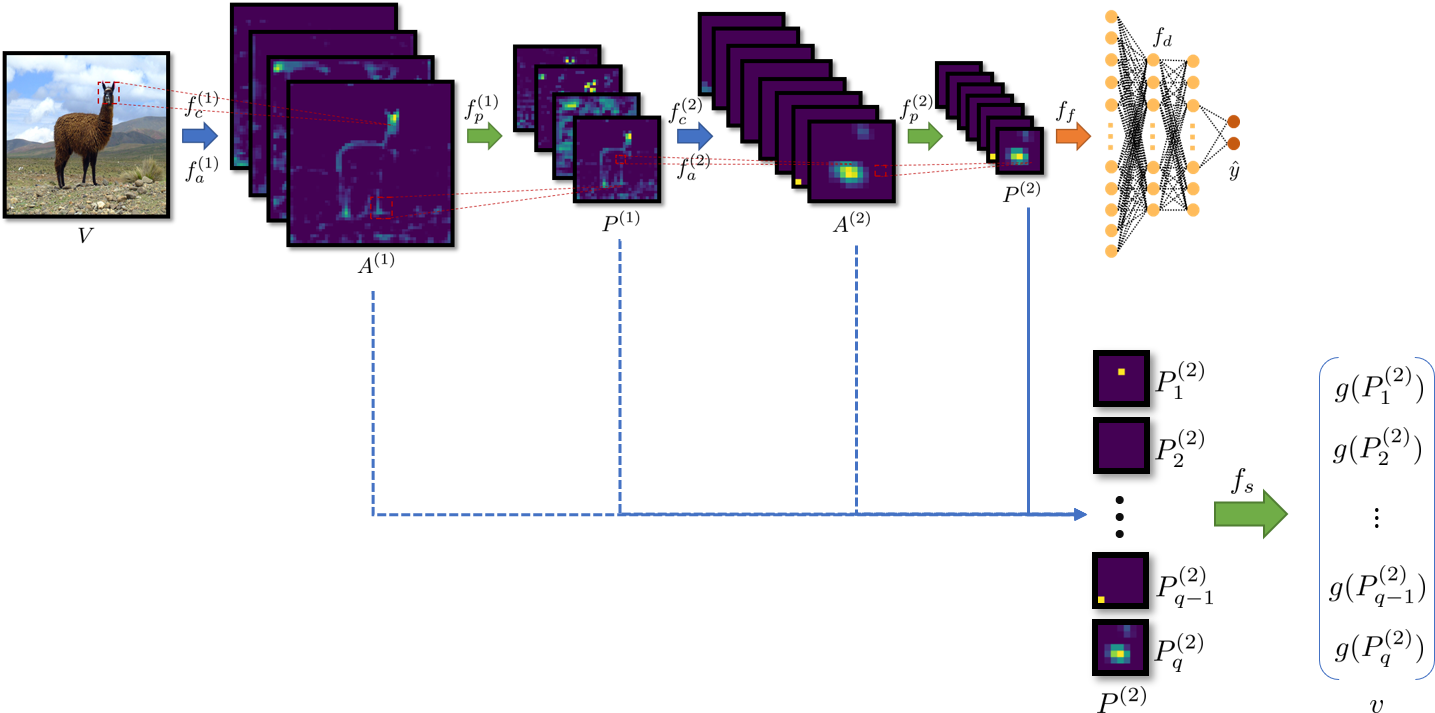}
    \caption{Illustration on how $f_s$ scalarizes each feature map matrix $P_j^{(2)}$ from a CNN sensor to derive a 1D vector of features $v$.}   
    \label{fig:feature_extract}
\end{figure}

Figure \ref{fig:feature_extract} illustrates our proposed feature extraction methodology and highlights key methodological flexibilities. First, we can choose from a variety of functions $g(\cdot)$ to scalarize each feature map matrix $P_j^{(\ell)}$ which follows from our discussion above. Second, we can select which feature map layer to extract from. In Section \ref{sec:cnn}, we observed how the initial convolutional block captures smaller length-scale visual patterns and how the deepest block highlights more sophisticated patterns/objects within an image. Thus, it may generally be advantageous to select the first block for deriving a feature space for conducting novelty detection on certain less sophisticated visual disturbances (e.g., blurring); similarly, the last block would be a natural choice for more complex disturbances. For a selected convolutional block, we can extract from the output of the pooling layer $P^{(\ell)}$ as shown in our above analysis, but we also have the methodological flexibility to instead use an intermediate output (i.e., $\Psi^{(\ell)}$ or $A^{(\ell)}$). For instance, we might select $\Psi^{(\ell)}$ if we want a feature space that is only influenced by the output of the convolutional filters (i.e., does not undergo the nonlinear and spatial reduction transformations induced by the activation and pooling layers).  

\subsection{Feature Refinement} \label{sec:feature_refine}

Once we have extracted the raw feature space following the methodology described in Section \ref{sec:feature_extract}, we will typically need to refine it to make it readily amendable for use with an OCC novelty detector. Hence, we will use a refinement mapping function $f_r : \mathbb{R}^q \mapsto \mathbb{R}^{d}$ that refines a feature vector $v \in \mathbb{R}^q$ into a refined variant $v' \in \mathbb{R}^d$ where $d \leq q$. Two common refinement avenues include feature scaling/normalization and dimension reduction. 

Feature scaling and normalization refers to a group of methods that linearly transform feature data such that each feature $v'_i$ is better conditioned. These transformations are typically carried out element-wise with a transformation function $s : \mathbb{R} \mapsto \mathbb{R}$ where common choices include:
\begin{equation*}
    \begin{gathered}
        s_\text{scale}(v_i) = \frac{v_i - v_{\min,i}}{v_{\max,i} - v_{\min,i}} \\
        s_\text{standard}(v_i) = \frac{v_i - v_{\mu,i}}{v_{\sigma,i}} \\
        s_\text{norm}(v_i) = \frac{v_i - v_{\mu,i}}{v_{\max,i} - v_{\min,i}}.
    \end{gathered}
\end{equation*}
The scaling transformation $s_\text{scale}(v_i)$ uses the minimum $v_{\min}$ and maximum $v_{\max}$ of the training data $\{v^{(k)} : k \in \mathcal{K}\}$ (computed element-wise) to scale each feature such that $v'_i \in [0, 1]$. Standardization uses $s_\text{standard}(v_i)$ to fit each feature to a standard Gaussian distribution $\mathcal{N}(0, 1)$ via the mean $v_{\mu,i}$ and standard deviation $v_{\sigma,i}$ of each feature. Finally, the normalization technique $s_\text{norm}(v_i)$ combines the scaling and standardization approaches to produce Gaussian-distributed features with a scaled standard deviation. Selecting an appropriate transformation function will depend on the nature of the training data. For instance, data with extreme outliers may be better suited for standardization so that the reminder of the training instances are not compressed into a small interval as would occur with scaling. 
\\

Moreover, we may need to reduce the size of $v$ if $q$ is large relative to the amount of available training instances and/or the desired OCC model. Following the discussion in Section \ref{sec:pca}, PCA is a natural choice since it employs a linear transformation that typically incurs a minor computational expense relative to the other components of our proposed framework. Here, we obtain a reduced feature vector $v' \in \mathbb{R}^d$ using the transformation $f_\text{pca}$ which is trained using the training data $\{v^{(k)} : k \in \mathcal{K}\}$. Note that the scaled/normalized data should be used. Selecting the appropriate amount $d$ of reduced features can be done by thresholding the total data variation retention against the number of principal components used (i.e., the value of $d$). Thus, we might refine our data via standardization and PCA:
\begin{equation}
    v' = f_r(v) = f_\text{pca}(s(v)).
\end{equation}
Alternative dimensionality reduction techniques include Sparse PCA, Linear Discriminant Analysis (LDA), t-Distributed Stochastic Neighbor Embedding (t-SNE), and Uniform Manifold Approximation and Projection (UMAP). We refer the reader to \cite{espadoto2019toward} for a comprehensive survey on reduction techniques.

\subsection{One-Class Classification} \label{sec:feature_occ}

We use the extracted and refined feature space to build an appropriate OCC novelty detector $f_\text{occ} : \mathbb{R}^d \mapsto \mathbb{R}$ which maps the feature vector $v' \in \mathbb{R}^d$ to the novelty prediction $\hat{h} \in \mathbb{R}$ which is thresholded to determine if $v'$ is novel. Our derived feature space follows a traditional 1D structure which can be readily used with any of the OCC techniques discussed in Section \ref{sec:occ}. For the remainder of this work, we will use OC-SVMs which have the advantage of incurring a relatively low computational cost to train and use; moreover, they are a standard choice for OCC. We refer the reader to \cite{khan2014one} to learn more about alternative OCC models. We train the OCC novelty detector using the refined feature data $\{v'^{(k)} : k \in \mathcal{K}\}$ where each $v'^{(k)}$ is derived from the corresponding training image $V^{(k)}$ that is used to train the CNN sensor. Here the training data is labeled as normal in its entirety. 
\\

\begin{figure}[!htb]
    \centering
    \includegraphics[width=0.8\textwidth]{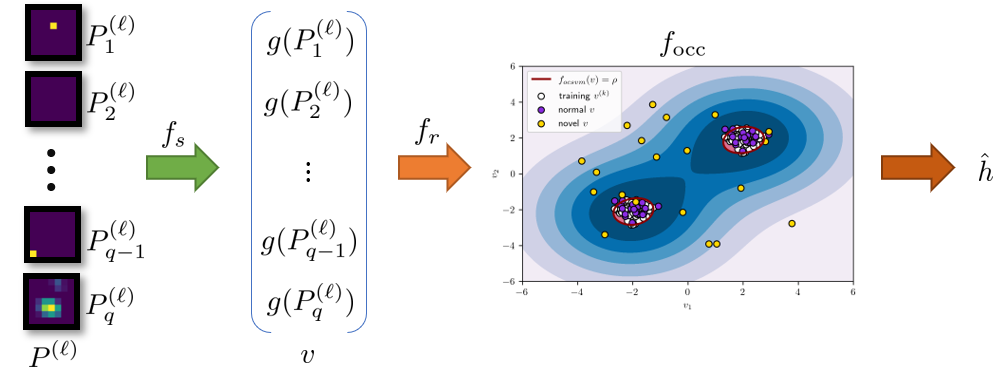}
    \caption{A summary of the SAFE-OCC novelty detection framework that operates on a feature map $P^{\ell}$ to produce a novelty signal $\hat{h}$.}
    \label{fig:feature_framework}
\end{figure}

Figure \ref{fig:feature_framework} summarizes our SAFE-OCC novelty detection framework that applies the feature extractor $f_s$ to a selected feature map $P^{(\ell)}$ to yield 1D features $v$ which are then transformed via $f_r$ to refined features $v'$ that are given to an OCC novelty detector $f_\text{occ}$ which gives the predicted novelty signal $\hat{h}$:
\begin{equation}
    \hat{h} = f_\text{occ}(f_r(f_s(P^{(\ell)}))).
\end{equation}
The novelty detection framework can be implemented in parallel using different configurations to augment the range of visual disturbances that can be detected. For instance, we can setup two independent systems that use the first and last feature map blocks, respectively. This would setup one novelty detector that is especially sensitive to lower length scale disturbances and another that is sensitive to more complex disturbances. We will demonstrate such a parallel system in the case study presented in Section \ref{sec:case1}. One potential extension of this concept would be in creating an ensemble of SAFE-OCC novelty detectors that employ a wide variety of architectures. Such an extension would likely expand the range of novelties that could be detected, and statistical information could be extracted from the results in similar manner to the neural network ensemble approach that was recently proposed in \cite{liu2021uncertainty} for uncertainty quantification. We leave the investigation of such an approach to future work.

\subsection{Incorporation into Control System}

\begin{figure}[!htb]
    \centering
    \includegraphics[width=0.8\textwidth]{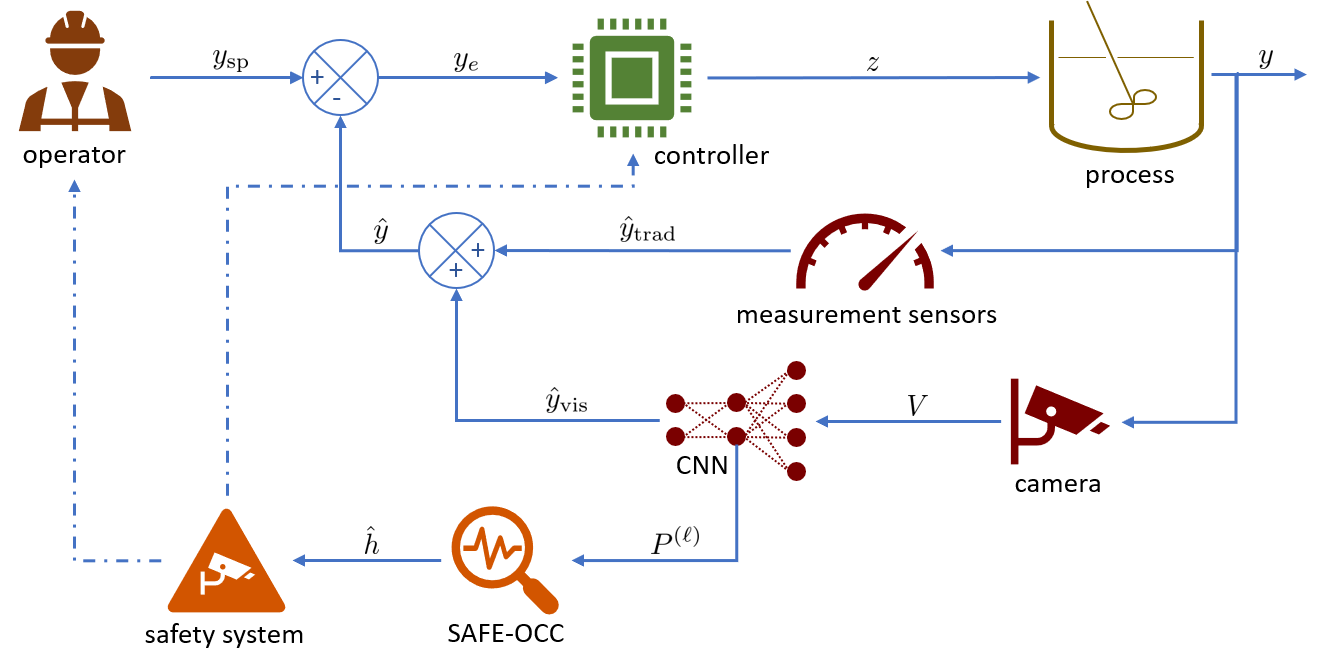}
    \caption{Control loop that incorporates the SAFE-OCC framework into a safety system that takes appropriate recourse action if novel input data is detected (making the CNN sensor unreliable).}
    \label{fig:safe_control}
\end{figure}

We now explain how the SAFE-OCC framework can be incorporated into the CNN sensor aided control architecture discussed in Section \ref{sec:compvis_control}. Figure \ref{fig:safe_control} illustrates this incorporation based on the control loop presented in Figure \ref{fig:vision_control}. Here, a feature map signal $P^{(\ell)}$ is extracted from the CNN sensor and feeds into the SAFE-OCC framework to obtain a novelty signal $\hat{h}$. This novelty signal is monitored by a safety system that uses pre-determined logic to invoke an appropriate response signal to the operator and/or the controller (e.g., raise an alarm and prompt the operator to revert to manual control). Such a safety system is crucial to mitigate the risk of injecting inaccurate CNN sensor measurements into a control system and incurring significant operation deviations. We implement this proposed control loop architecture in the case study presented in Section \ref{sec:case2}.

\section{Case Studies} \label{sec:examples}

We now present illustrative case studies to highlight the effectiveness of the SAFE-OCC novelty detection framework and to demonstrate its utility for control systems that incorporate CNNs to measure process states from image data. We use simulated environments from \texttt{OpenAI-Gym} \cite{1606.01540} to serve as the processes we seek to control. Moreover, we implement our CNN models using \texttt{TensorFlow} and we use \texttt{SciKit-Learn} to implement OC-SVMs and PCA reduction. Finally, we augment the training data (collected using \texttt{OpenAI-Gym}) and simulate disturbances using \texttt{ImgAug}. 

\subsection{Simple SAFE-OCC Novelty Detection} \label{sec:case1}

In this case study, we apply the SAFE-OCC novelty detection framework to a simple process system and investigate its effectiveness over a range of CNN sensors and visual disturbances. 

\subsubsection{CNN Sensor Training}

We use the \texttt{Pendulum-v0} environment from \texttt{OpenAI-Gym} to generate images sets $\{V^{(k)}: k \in \mathcal{K}\}$ that map to state variables $\{y^{(k)} : k \in \mathcal{K}\}$. This environment corresponds to the classical inverted pendulum swing-up problem where a pendulum with a fixed axis of rotation can be swung left or right via a controllable torque input. We generate our labeled image data via 200 simulations that employ random control inputs to produce 7,750 images that are divided in a 70:20:10 split for training, validation, and test data sets, respectively. These images are extracted from the frames of each simulation video. Moreover, each image corresponds to the state variables:
\begin{equation}
    y = \begin{bmatrix} \sin(\theta) \\ \cos(\theta) \end{bmatrix} \in \mathbb{R}^2
\end{equation}
where $\theta$ is the angle of the pendulum relative to an upright vertical position.

\begin{table}[!htb]
    \centering
    \caption{The \texttt{ImgAug} image augmentation types used to simulate disturbances.}
    \begin{tabular}{| c | c |} 
        \hline
        \texttt{ImgAug} Type & Disturbance Description \\
        \hline
        \texttt{Cutout} & Random blockages \\
        \texttt{DefocusBlur} & Out-of-focus blurring \\
        \texttt{Fog} & Simulated fog \\
        \texttt{GaussianNoise} & Gaussian noise \\
        \texttt{PerspectiveTransform} & Shift in camera direction \\
        \texttt{Spatter} & Lens splattering \\
        \hline
    \end{tabular}
    \label{tab:disturbance_types}
\end{table}

\begin{figure}[!htb]
     \centering
     \begin{subfigure}[b]{0.15\textwidth}
         \centering
         \includegraphics[width=\textwidth]{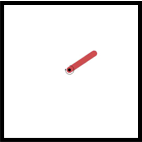}
         \caption{Original }
     \end{subfigure}
     \quad
     \begin{subfigure}[b]{0.15\textwidth}
         \centering
         \includegraphics[width=\textwidth]{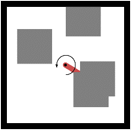}
         \caption{Blockages }
     \end{subfigure}
     \quad
     \begin{subfigure}[b]{0.15\textwidth}
         \centering
         \includegraphics[width=\textwidth]{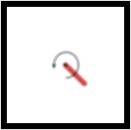}
         \caption{Blurred }
     \end{subfigure}
     \quad
     \begin{subfigure}[b]{0.15\textwidth}
         \centering
         \includegraphics[width=\textwidth]{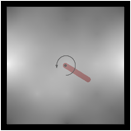}
         \caption{Fogged }
     \end{subfigure}
     \\
     \begin{subfigure}[b]{0.15\textwidth}
         \centering
         \includegraphics[width=\textwidth]{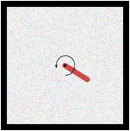}
         \caption{Noised }
     \end{subfigure}
     \quad
     \begin{subfigure}[b]{0.15\textwidth}
         \centering
         \includegraphics[width=\textwidth]{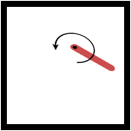}
         \caption{Shifted }
     \end{subfigure}
     \quad
     \begin{subfigure}[b]{0.15\textwidth}
         \centering
         \includegraphics[width=\textwidth]{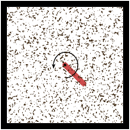}
         \caption{Splattered }
     \end{subfigure}
    \caption{Samples of the disturbances highlighted in Table \ref{tab:disturbance_types}. }
    \label{fig:DistLabel}
\end{figure}

We augment our image data using \texttt{ImgAug} to simulate a variety of visual disturbances. In particular, we use the disturbance types described in Table \ref{tab:disturbance_types}. Figure \ref{fig:DistLabel} exemplifies the aforementioned disturbances that augment the dataset. These are intended to represent a typical range of visual disturbances that might transpire in process control applications. 

\begin{figure}[!htb]
    \centering
    \includegraphics[width=0.85\textwidth]{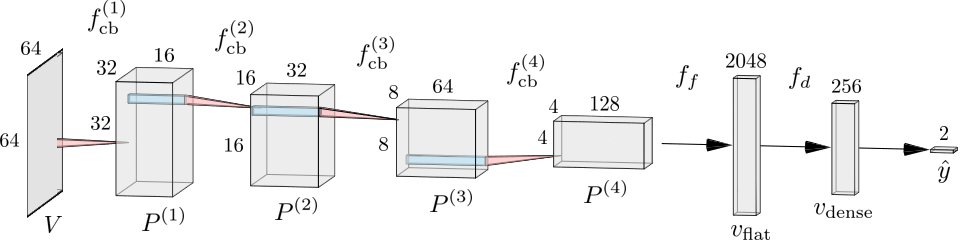}
    \caption{Schematic of the CNN used to predict the states $\hat{y}$ given a gray-scale image $V$ in the case study presented in Section \ref{sec:case1}.}   
    \label{fig:case1_cnn}
\end{figure}

The CNN architecture we use for this study is summarized in Figure \ref{fig:case1_cnn}. We employ four blocks $f_\text{cb}^{\ell}, \ \ell \in \{1, 4\},$, that each consist of a single convolution, activation, and pooling layer. Here, the dimension of each block reduction are provided in Figure \eqref{fig:case1_cnn}, and we use ReLU activation functions and $2 \times 2$ max-pooling operators. All convolutional operators use $n_u = 3$ and zero-padding is used. We flatten the final feature map output $P^{(4)}$ via $f_f$ and map it through the dense layer $f_d$ to yield the predicted state variable $\hat{y}$. We implement this CNN sensor architecture in \texttt{TensorFlow} using the Adam optimizer \cite{kingma2014adam} to train each sensor. The learning rate hyper-parameter is chosen for each sensor via an enumeration study juxtaposing the SSE training loss after one iteration. Each sensor is training until the training loss is observed graphically to level-off, and we ensure that it does surpass the validation loss (to prevent over-fitting). In total, we train seven CNN sensors where each one is given a different training image set. Sensor A is training using only the original simulation images, while Sensors B-G are trained using the blockage, blurred, fogged, noised, shifted, and spattered images, respectively, in addition to using the original images. 

\subsubsection{CNN Sensor Accuracy Assessment}

\begin{table}[!htb]
    \centering
    \caption{Mean $L_{2}$ prediction error for each CNN sensor relative to each image test set. The shaded and the non-shaded values correspond to normal and novel test data, respectively (relative to the training data of each sensor).}
    \begin{tabular}{| c || c | c | c | c | c | c | c |} 
        \hline
        \multirow{2}{*}{Sensor} & \multicolumn{7}{c|}{Test Data}  \\
        \cline{2-8}
        & Original & Blockages & Blurred & Fogged & Noised & Shifted & Splattered \\
        \hline \cline{2-8}
        A & \cellcolor{black!15}0.024 & 1.092 & 0.165 & 1.141 & 0.598 & 0.750 & 0.890 \\ 
        B & \cellcolor{black!15}0.060 & \cellcolor{black!15}0.293 & 0.197 & 0.620 & 0.148 & 0.958 & 0.903 \\ 
        C & \cellcolor{black!15}0.023 & 1.507 & \cellcolor{black!15}0.021 & 1.493 & 0.540 & 0.760 & 1.047 \\ 
        D & \cellcolor{black!15}0.029 & 1.094 & 0.143 & \cellcolor{black!15}0.066 & 0.382 & 0.972 & 0.576 \\ 
        E & \cellcolor{black!15}0.028 & 1.383 & 0.161 & 1.022 & \cellcolor{black!15}0.031 & 0.774 & 0.860 \\ 
        F & \cellcolor{black!15}0.046 & 1.209 & 0.138 & 1.483 & 0.707 & \cellcolor{black!15}0.199 & 0.703 \\ 
        G & \cellcolor{black!15}0.019 & 1.296 & 0.069 & 0.975 & 0.155 & 0.660 & \cellcolor{black!15}0.055 \\
        \hline
    \end{tabular}
    \label{tab:MeanError}
\end{table}

We subject the trained CNN sensors to each disturbance test data set and the mean $L_{2}$ predictions errors ($|\hat{y} - y|_2$) are recorded in Table \ref{tab:MeanError}. We observe that the normal image sets for each sensor (i.e., test images that are similar to the training data) generally incur low prediction errors, while the novel images tend to induce significant prediction error. A couple of apparent outliers to this trend are the mean blockage and shifted errors for models B and F, respectively. This occurs since certain images in these data sets are highly adversarial (e.g., certain blockages completely cover the pendulum) making accurate predictions impossible, but we note that mean predictions are significantly lower than those incurred by their counterpart models. These results illustrate how image augmentation is typically effective in mitigating the effects of disturbances, but novel disturbances (which correspond to extrapolating predictions) tend to induce large prediction errors that can have significant ramifications for process control applications as discussed above.

\begin{figure}[!htb]
    \centering
    \includegraphics[width=0.5\textwidth]{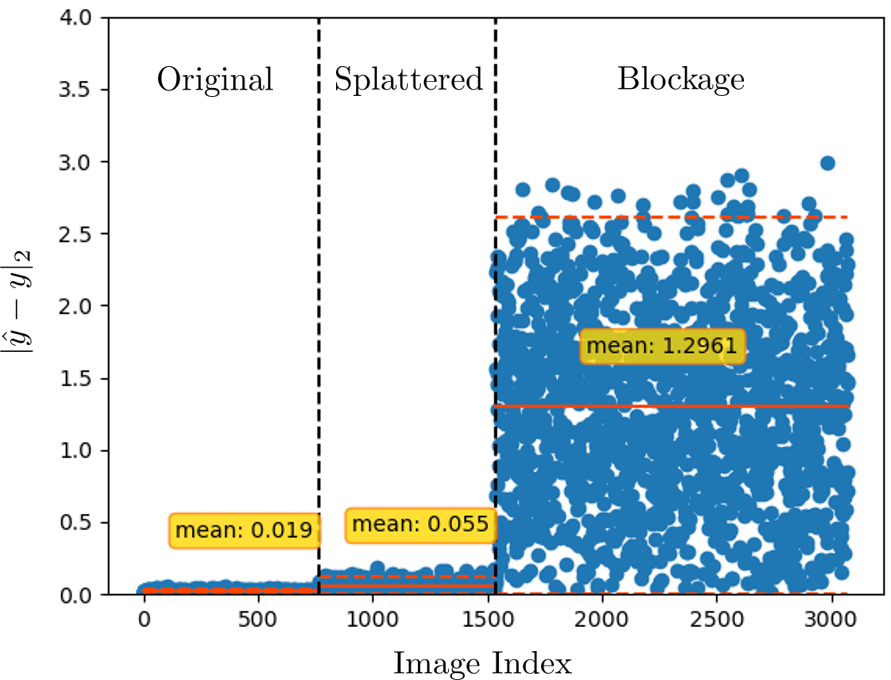}
    \caption{Depiction of prediction errors incurred by subjecting CNN Sensor G to original and splattered (normal) images vs. blockage (novel) images. The novel images induce significantly more error on average.}   
    \label{fig:splatblock}
\end{figure}

Figure \ref{fig:splatblock} illustrates this principle further by showing the individual prediction errors observed by subjecting Sensor G (trained with original and splattered images) to original, splattered, and blockage image test data sets. Here, we readily observe how each original and splattered image incurs a significantly lower prediction error relative to the blockage images. It is also apparent that the novel blockage images exhibit a significantly higher mean and variance in their prediction errors. These can be interpreted as situations in which the sensor makes uninformed guesses for the state values which leads to a large spread in the individual prediction accuracy. 

\subsubsection{SAFE-OCC Configuration Study}

The SAFE-OCC novelty detection framework provides flexibility in how it can be configured for a particular CNN sensor of interest. We investigate this flexibility in selecting the input feature map $P^{(\ell)}$, the scalarization function $g(\cdot)$, and the feature refinement transformation function $s(\cdot)$. We exhaustively explore the combinatorics of these design choices in conjunction with CNN Sensor A and subject each one to the disturbance types to assess the effectiveness of the framework to classify normal and novel images. Our investigation considers the convolutional, activation, and max-pooling outputs of the first and last convolutional blocks. Moreover, we scalarize each of the feature map outputs using either $g_\text{2d2pca}$ or $g_\text{max}$. Finally, we explore leaving the resulting feature vectors unrefined, refined using $s_\text{standard}$, or refined using $s_\text{scale}$. Each resulting configuration is then paired with a OC-SVM novelty detector which is trained via the training data of CNN Sensor A. Each OC-SVM (implemented in \texttt{SciKit-Learn}) uses a Gaussian kernel function with $\gamma = 1/n_v$ and sets the hyper-parameter $\nu = 0.0001$ since we are confident that the training data does not contain mislabeled normal images (all other hyperparameters are kept at their defaults). We then subject each configuration to each test image set, and the results are provided in Table \ref{tab:CleanImageStudy} (Appendix \ref{sec:appendix}). 
\\

Table \ref{tab:CleanImageStudy} reveals several trends regarding the performance of certain SAFE-OCC configurations relative to this case study. We find that the $g_\text{2d2pca}$ scalarization function performs better in combination with a feature map output from the last convolutional block, and that $g_\text{max}$ performs well in combination with the first layer outputs. Moreover, we observe that the refinement transformation $s_\text{standard}$ performs better than $s_\text{scale}$ across all the configurations considered in this study, and that the convolutional and pooling feature map outputs outperform their activation layer counterparts. Another key observation is that the configurations that use feature maps from the first convolutional block tend to better detect novel blurring, while the remaining disturbance types are detected more readily in general by configurations that use feature maps from the last convolutional block. This is consistent with the discussion in Section \ref{sec:cnn} on how deeper convolutional blocks extract more specialized features over larger lengthscales. This highlights the utility of using diverse SAFE-OCC configurations to detect a range of anomaly types.

\begin{figure}[!htb]
     \centering
     \begin{subfigure}[b]{0.31\textwidth}
         \centering
         \includegraphics[width=\textwidth]{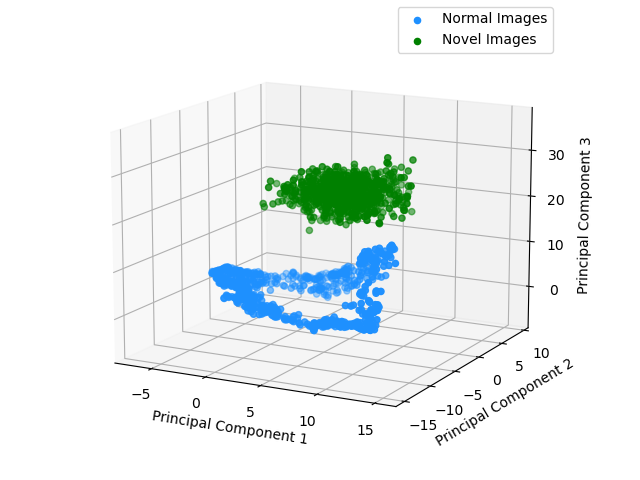}
         \caption{Original and Splattered}
     \end{subfigure}
     \quad
     \begin{subfigure}[b]{0.31\textwidth}
         \centering
         \includegraphics[width=\textwidth]{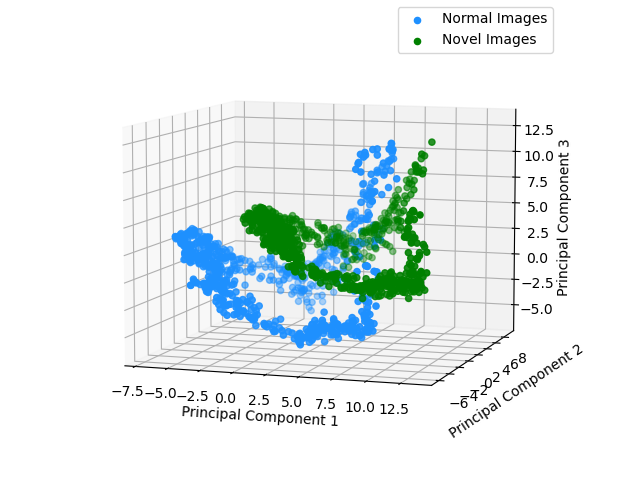}
         \caption{Original and Noised}
     \end{subfigure}
     \quad
     \begin{subfigure}[b]{0.31\textwidth}
         \centering
         \includegraphics[width=\textwidth]{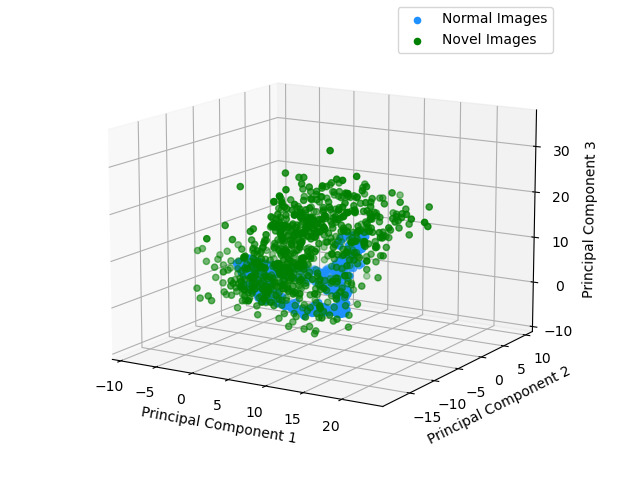}
         \caption{Original and Blocked}
     \end{subfigure}
    \caption{Representative collection of 3D PCA projections of the feature spaces studied in Table \ref{tab:CleanImageStudy}. These correspond to the configuration that uses $p^{(4)}$, $g_\text{2d2pca}$, and $s_\text{standard}$. The feature points corresponding to original (normal) and disturbed (novel) images are near perfectly separable.}
    \label{fig:3DPCAplots}
\end{figure}

Figure \ref{fig:3DPCAplots} illustrates the feature spaces (projected into a reduced PCA space) of the SAFE-OCC configuration that scalarizes the last pooling layer feature map with $g_\text{2d2pca}$ and standardizes the features via $s_\text{standard}$. Interestingly, we observe that a clear separation between the novel and normal data points in this greatly reduced feature space. This adds some visual intuition as to why the full SAFE-OCC derived feature spaces are generally effective for novelty detection. Moreover, this highlights that using dimension reduction techniques detailed in Section \ref{sec:feature_refine} such as PCA in combination with the SAFE-OCC framework can derive high-fidelity feature spaces that use a low number of features. This property can help reduce the computational burden of novelty detection and enable the use of more sophisticated novelty detectors.
\\

\begin{figure}[!htb]
    \centering
    \includegraphics[width=\textwidth]{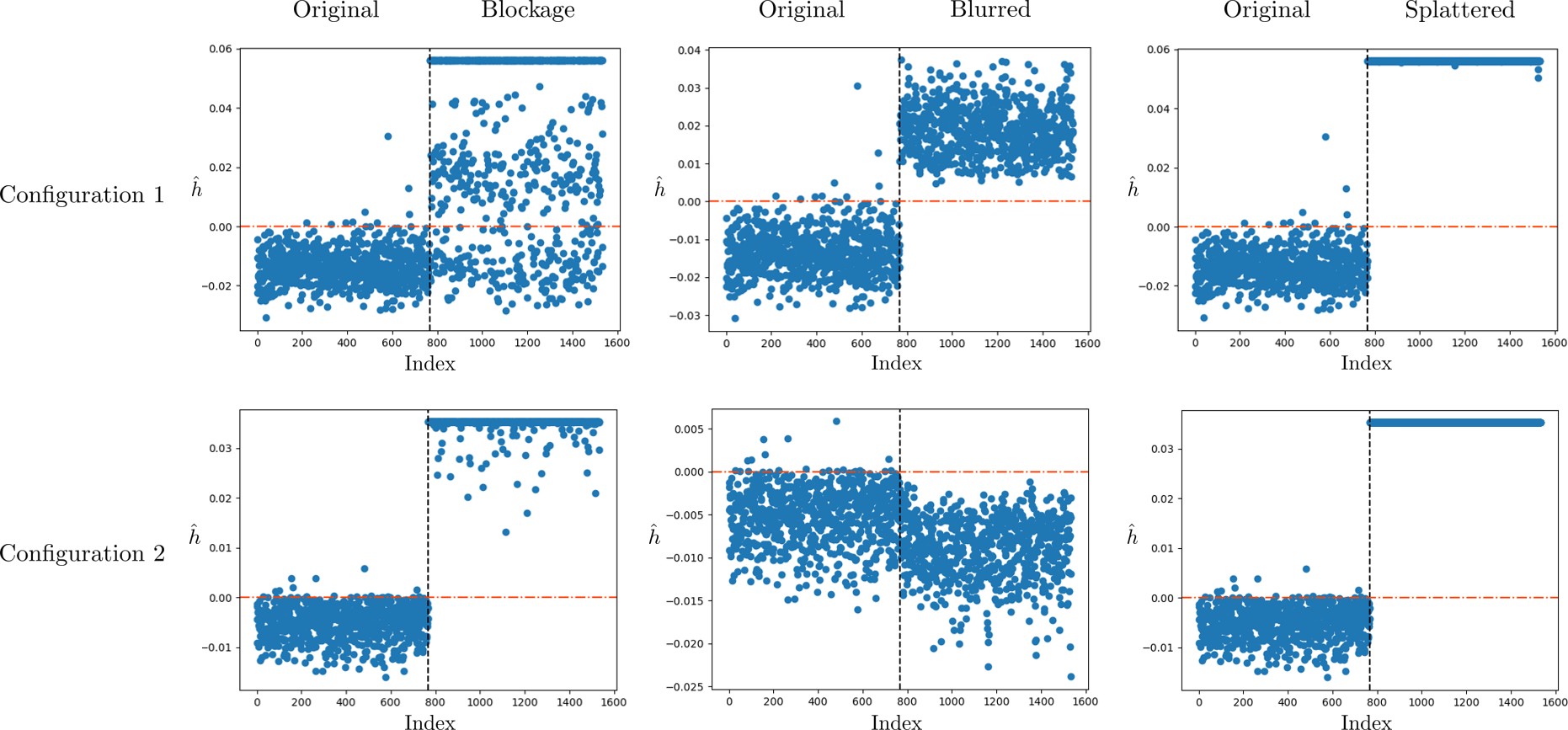}
    \caption{Novelty predictions $\hat{h}$ of SAFE-OCC Configurations 1 and 2 applied to CNN Sensor A in conjunction with the original (the training data), blockage, blurred, and splattered datasets. Both configurations compensate for the deficiencies of each other.}
    \label{fig:NoveltyPlots}
\end{figure}

With the results Table \ref{tab:CleanImageStudy}, we select two SAFE-OCC configurations that we will apply to each CNN sensor in Section \ref{sec:case1_safe_occ} of this case study. One configuration applies the scalarization function $g_\text{max}$ to the max-pooling layer output $P^{(1)}$ and refines it element-wise via $s_\text{standard}$; we refer to it as Configuration 1. Configuration 2 employs the max-pooling output $P^{(4)}$ with $g_\text{2d2pca}$ and $s_\text{standard}$. Both of the feature spaces derived from these configurations are used in combination with a OC-SVM novelty detector that uses a Gaussian kernel. Figure \ref{fig:NoveltyPlots} shows the predicted novelty output $\hat{h}$ (instances above the $\rho=0$ threshold are classified as novel) of each SAFE-OCC configuration implemented on CNN Sensor A with respect to the original, blockage, blurred, and splattering test datasets. We readily observe the complementary nature of the two configurations in how Configuration 1 is able to effectively identify the blurred and splattered disturbances and how Configuration 2 is classifies the blockage and splattered disturbances to high accuracy. In this regard, Configurations 1 and 2 complement and compensate for the deficiencies of each other.

Furthermore, Figure \ref{fig:NoveltyPlots} visually motivates a practical consideration for novelty detection: shifting the novelty detection threshold $\rho$ by a small tolerance $\epsilon \in \mathbb{R}$. This helps to decrease the frequency of false positives (erroneously labeling normal images as novel) which can occur with images that immediately border the learned OC-SVM boundary. In this case, we can see that raising the thresholds in Figure \ref{fig:NoveltyPlots} from 0 to 0.1 would increase the normal accuracy by having nearly all the normal original images correspond novelty predictions below the threshold. Naturally, the threshold can also be adjusted by more finely tuning the OCC novelty detector.

\subsubsection{SAFE-OCC Novelty Detection} \label{sec:case1_safe_occ}
As discussed in Section \ref{sec:feature_occ}, multiple SAFE-OCC framework configurations can be implemented in parallel to detect a larger range of disturbances. The choice of our two configurations will further illustrate this as we apply them to all seven of the CNN sensors featured in this case study.
\\

Table \ref{tab:FirstPool} shows the classification accuracies achieved by applying Configuration 1 to all the CNN sensors and subjecting the resulting SAFE-OCC framework to all the test datasets. Note that when the test datasets that relate to the data used to train a particular sensor, the normal accuracy is considered; whereas the other datasets should be identified as novel and their novel accuracy is shown. Configuration 1 is able to achieve high novelty accuracy for the lower length-scale disturbance types (i.e., blurred, fogged, noised, and splattered). The blurred novelty accuracy for Sensor B provides one exception, but this can likely be attributed to the training data containing some highly adversarial blockage images (i.e., ones where the pendulum is completely obscured) and in practice these should be excluded from the training set. Importantly, the normal accuracy is high in all cases (meaning that false-positives are infrequent). Moreover, in other tests we found that normal accuracy can be made near perfect by shifting the novelty threshold by a small tolerance as discussed above. 

\begin{table}[!htb]
    \centering
    \caption{The classification accuracies (in \%) of SAFE-OCC Configuration 1 for each CNN sensor. The shaded and non-shaded values denote the normal and novelties accuracies, respectively, (with respect to the training data of each sensor).}
    \begin{tabular}{|c||c|c|c|c|c|c|c|}
    \hline
        \multirow{2}{*}{Sensor} & \multicolumn{7}{c|}{Test Data}  \\
        \cline{2-8}
          & Original & Blockages & Blurred & Fogged & Noised & Shifted & Splattered  \\  
        \hline \cline{2-8}
        A & \cellcolor{black!15}98.57 & 74.61 & 100.00 & 100.00 & 100.00 & 74.87 & 100.00 \\
        B & \cellcolor{black!15}95.70 & \cellcolor{black!15}95.70 & 21.68 & 100.00 & 100.00 & 11.65 & 100.00 \\ 
        C & \cellcolor{black!15}98.76 & 82.49 & \cellcolor{black!15}98.76 & 100.00 & 100.00 & 44.92 & 100.00 \\ 
        D & \cellcolor{black!15}99.74 & 67.12 & 87.89 & \cellcolor{black!15}99.74 & 1.30 & 27.93 & 100.00 \\ 
        E & \cellcolor{black!15}99.35 & 63.15 & 100.00 & 100.00 & \cellcolor{black!15}99.35 & 28.58 & 100.00 \\ 
        F & \cellcolor{black!15}97.20 & 68.68 & 97.85 & 100.00 & 100.00 & \cellcolor{black!15}97.20 & 100.00 \\ 
        G & \cellcolor{black!15}98.44 & 48.50 & 100.00 & 100.00 & 100.00 & 26.17 & \cellcolor{black!15}98.44 \\
        \hline
    \end{tabular}
    \label{tab:FirstPool}
\end{table}

Table \ref{tab:LastPool} shows the classification accuracies that come from applying Configuration 2 to all the CNN sensors. As expected, we observe that is performs better than Configuration 2 on the longer length-scale disturbance types (i.e., blockage and shifted). Again, high normal accuracies are obtained which can be made near perfect by adjusting the novelty threshold by a small numerical tolerance. This exemplifies other empirical studies we have done which show the SAFE-OCC novelty detection framework is not particularly prone to false-positives. This property is what makes the SAFE-OCC framework readily parallelizable where novel detections can be summarized by simply employing the union of the results. 

\begin{table}[!htb]
    \centering
    \caption{The classification accuracies (in \%) of SAFE-OCC Configuration 2 for each CNN sensor. The shaded and non-shaded values denote the normal and novelty accuracies, respectively, (with respect to the training data of each sensor).}.
    \begin{tabular}{|c||c|c|c|c|c|c|c|}
    \hline
        \multirow{2}{*}{Sensor} & \multicolumn{7}{c|}{Test Data}  \\
        \cline{2-8}
          & Original & Blockages & Blurred & Fogged & Noised & Shifted & Splattered  \\  
        \hline \cline{2-8}
        A & \cellcolor{black!15}97.27 & 100.00 & 0.00 & 100.00 & 100.00 & 96.35 & 100.00 \\
        B & \cellcolor{black!15}97.07 & \cellcolor{black!15}97.07 & 2.02 & 100.00 & 17.84 & 38.09 & 100.00 \\
        C & \cellcolor{black!15}98.31 & 100.00 & \cellcolor{black!15}98.31 & 100.00 & 100.00 & 98.05 & 100.00 \\ 
        D & \cellcolor{black!15}99.02 & 99.93 & 0.13 & \cellcolor{black!15}99.02 & 39.06 & 74.48 & 100.00 \\ 
        E & \cellcolor{black!15}98.44 & 100.00 & 3.06 & 100.00 & \cellcolor{black!15}98.44 & 91.67 & 100.00 \\ 
        F & \cellcolor{black!15}96.29 & 98.31 & 0.00 & 100.00 & 100.00 & \cellcolor{black!15}96.29 & 100.00 \\ 
        G & \cellcolor{black!15}100.00 & 99.87 & 10.29 & 100.00 & 100.00 & 42.45 & \cellcolor{black!15}100.00  \\
        \hline
    \end{tabular}
    \label{tab:LastPool}
\end{table}

In this case study, the integration of these configurations yields a novelty detection system that is more effective than what could be achieved by using a single configuration. From Tables \ref{tab:FirstPool} and \ref{tab:LastPool} we can see that Configuration 1 struggles with detecting shifted and blockage disturbances while Configuration 2 struggles with blurred disturbances. Thus, the two respective SAFE-OCC configurations are complementary to each other. The resulting parallel framework still exhibits lower novelty accuracies for certain shifted disturbances. This in part can be attributed to some shifted images being very minorly disturbed (and thus weakly novel). Moreover, we can envision how adding more parallel configurations would increase the effectiveness of the SAFE-OCC novelty detection system.
\\

We have demonstrated with multiple CNN sensors that the SAFE-OCC framework can be quite effective at identifying novel image data that typically induces high sensor prediction errors. Moreover, we have illustrated how the tendency of the SAFE-OCC framework to rarely mislabel normal images in combination with its highly flexible nature, makes it readily parallelizable to produce a novelty detection system that is highly effective at identifying novel image data relative to a CNN sensor of interest.

\subsection{SAFE-OCC Aided Cart-Pole Control} \label{sec:case2}

We apply the SAFE-OCC aided control loop featured in Figure \ref{fig:safe_control} to control the \texttt{CartPole-v1} environment from \texttt{OpenAI-Gym} which corresponds to the classic cart-pole control problem introduced in \cite{barto1983neuronlike}. Here, we seek to balance a pendulum above a cart which we can move either right or left at a mixed rate. For simplicity in example, we consider the angle of the pendulum (measured in degrees relative to vertical alignment) as the state variable $y \in [-180, 180] \subset \mathbb{R}$ (ignoring the position of the cart), and we take the cart movement direction to be the control variable $z \in \{0, 1\} \subset \mathbb{Z}$ (where 0 is left and 1 is right). Thus, we have a single-input single-output (SISO) process to control. With this simplification, we implement a PID controller with a derivative filter whose control output is mapped through the sigmoid function and rounded to yield the binary control variable $z$. Note that such a controller is not robust for this difficult process, but it is able to maintain adequate control for the conditions considered in this case study.
\\

\begin{figure}[!htb]
    \centering
    \includegraphics[width=0.9\textwidth]{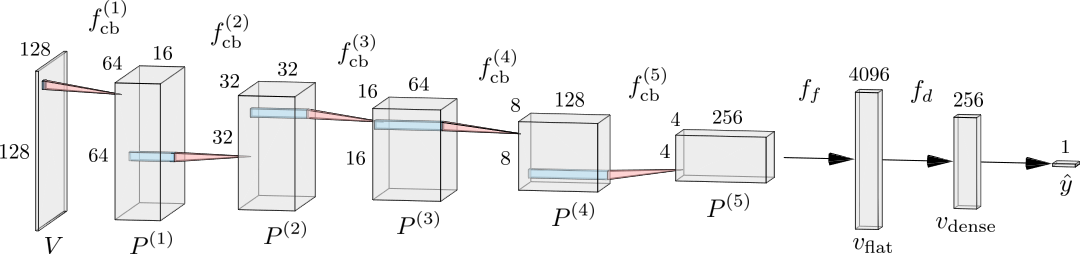}
    \caption{A schematic of the CNN sensor used to predict the pendulum angle $\hat{y}$ given a gray-scale simulation image $V$ in the case study presented in Section \ref{sec:case2}.}   
    \label{fig:case2_cnn}
\end{figure}

The predicted angle of pendulum $\hat{y}$ is estimated using a CNN sensor that takes gray-scale images $V \in \mathbb{R}^{128 \times 128}$ from the \texttt{OpenAI-Gym} simulation window as input. The accuracy of these predictions is readily assessed using the true value of $y$ which is simultaneously provided during the simulation. The structure of our CNN sensor is summarized in Figure \ref{fig:case2_cnn}. It uses five convolutional blocks $f_\text{cb}^{\ell}, \ \ell \in \{1, 5\},$ and in all other aspects is setup and trained in like manner to the CNN sensors in Section \ref{sec:case1}. We generate the training data via 1,000 simulations that use uniform random control input and terminate when the pendulum angle/position surpasses its default limit or after 200 time-steps. This produces 22,770 labeled images which we augment with the \texttt{Fog} disturbance from \texttt{ImgAug} to yield a training data set with 45,540 labeled images. This data is split randomly into 70:20:10 portions that correspond to training, validation, and test datasets, respectively.

The SAFE-OCC novelty detection framework is implemented using $P^{(5)} \in \mathbb{R}^{4 \times 4 \times 256}$ where each feature map $P_j^{(5)} \in \mathbb{R}^{4 \times 4}$ is scalarized using $g_\text{2d2pca}$. The resulting feature vector $v \in \mathbb{R}^{256}$ is then standardized element-wise using $s_\text{standard}$ to yield the refined feature vector $v' \in \mathbb{R}^{256}$. We implement a OC-SVM in \texttt{SciKit-Learn} that uses a radial basis kernel to learn a boundary around the training data. The framework is trained using the same set of training images (including the augmented ones) that were used to train the CNN sensor.  

\begin{figure}[!htb]
     \centering
     \begin{subfigure}[b]{0.2\textwidth}
         \centering
         \includegraphics[width=\textwidth]{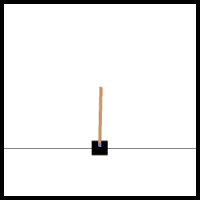}
         \caption{Simulation 1}
     \end{subfigure}
     \quad
     \begin{subfigure}[b]{0.2\textwidth}
         \centering
         \includegraphics[width=\textwidth]{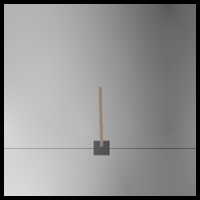}
         \caption{Simulation 2}
     \end{subfigure}
     \quad
     \begin{subfigure}[b]{0.2\textwidth}
         \centering
         \includegraphics[width=\textwidth]{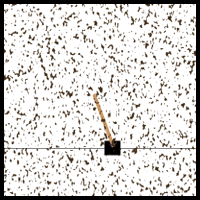}
         \caption{Simulation 3}
     \end{subfigure}
     \quad
     \begin{subfigure}[b]{0.2\textwidth}
         \centering
         \includegraphics[width=\textwidth]{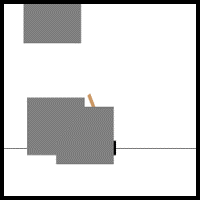}
         \caption{Simulation 4}
     \end{subfigure}
    \caption{Representative snapshots from the four simulations used in the cart-pole case study.}
    \label{fig:cart_images}
\end{figure}

We conduct four simulations: a base case that uses unperturbed images and three others that invoke a particular simulated visual disturbance after 150 time-steps. The three disturbance types are produced via \texttt{ImgAug} using the \texttt{Fog}, \texttt{Spatter}, and \texttt{Cutout} methods which correspond to fog, splattering, and square blockages, respectively. Figure \ref{fig:cart_images} shows representative images of these simulations. Each simulation is allotted a maximum duration of 400 time-steps or terminates prematurely when the pendulum rotates more than $180^\circ$ in either direction (it falls directly below the cart). Moreover, the setpoint $y_\text{sp}$ of the PID controller is set to $0^\circ$ (having the pendulum perfectly vertical).

\begin{figure}[!htb]
     \centering
     \begin{subfigure}[b]{0.45\textwidth}
         \centering
         \includegraphics[width=\textwidth]{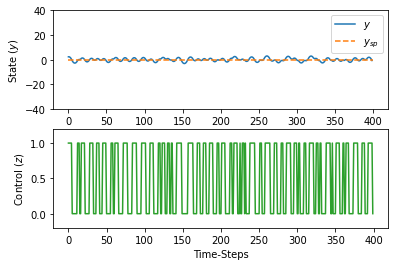}
         \caption{Simulation 1 Control Response}
         \label{fig:sim1_control}
     \end{subfigure}
     \quad
     \begin{subfigure}[b]{0.45\textwidth}
         \centering
         \includegraphics[width=\textwidth]{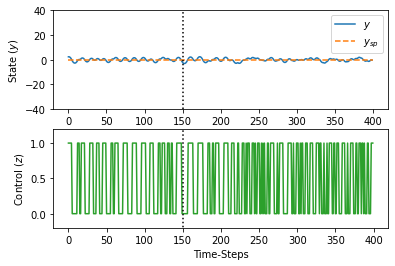}
         \caption{Simulation 2 Control Response}
         \label{fig:sim2_control}
     \end{subfigure}
     \\
     \begin{subfigure}[b]{0.45\textwidth}
         \centering
         \includegraphics[width=\textwidth]{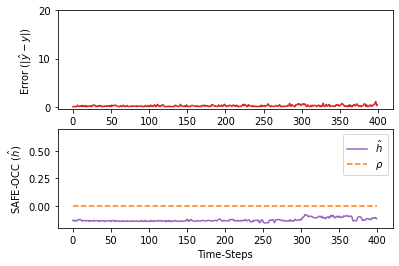}
         \caption{Simulation 1 Novelty Response}
         \label{fig:sim1_error}
     \end{subfigure}
     \quad
     \begin{subfigure}[b]{0.45\textwidth}
         \centering
         \includegraphics[width=\textwidth]{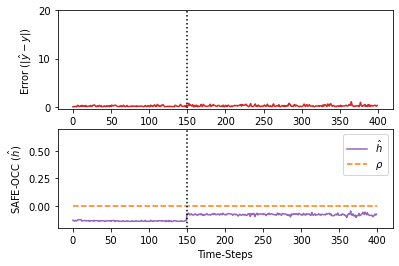}
         \caption{Simulation 2 Novelty Response}
         \label{fig:sim2_error}
     \end{subfigure}
    \caption{The control response trajectories of Simulations 1 and 2. The vertical dotted line at time-step 150 indicates when the fog disturbance is first introduced to Simulation 2. Effective control is maintained with these normal images as supported by the SAFE-OCC framework.}
    \label{fig:cart_control_normal}
\end{figure}

Figure \ref{fig:cart_control_normal} shows the responses exhibited in Simulations 1 and 2. In Figures \ref{fig:sim1_control} and \ref{fig:sim2_control}, we observe that effective control in tracking the setpoint is achieved for both simulations. This behavior can be attributed to the CNN sensor being trained with clear and fogged images which means that its predictions $\hat{y}$ incur a low error relative to $y$ as shown in Figures \ref{fig:sim1_error} and \ref{fig:sim2_error}. Moreover, the predicted novelty score $\hat{h}$ output by the SAFE-OCC framework is consistent with this observation, since its response trajectory remains below the novelty threshold $\rho$ in both simulations which means it correctly identified the incoming images as normal. Note that, in a real-world control system we cannot assess the error $|\hat{y} - y|$ in real-time with the true state $y$ being unknown. This is why it is essential to have an effective novelty detector that correlates with novel situations relative to the CNN sensor that induce high prediction error.
\\

\begin{figure}[!htb]
    \centering
    \begin{subfigure}[b]{0.45\textwidth}
        \centering
        \includegraphics[width=\textwidth]{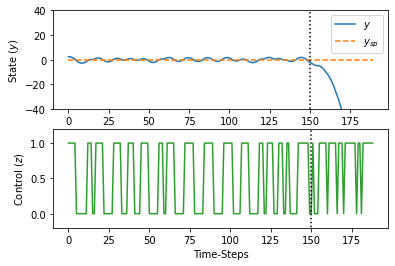}
        \caption{Simulation 3 Control Response}
        \label{fig:sim3_control}
    \end{subfigure}
    \quad
    \begin{subfigure}[b]{0.45\textwidth}
        \centering
        \includegraphics[width=\textwidth]{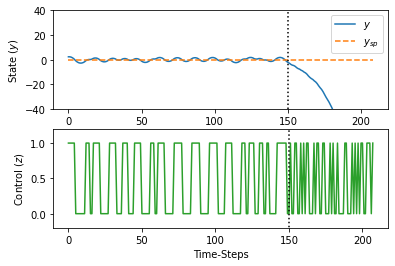}
        \caption{Simulation 4 Control Response}
        \label{fig:sim4_control}
    \end{subfigure}
    \\
    \begin{subfigure}[b]{0.45\textwidth}
        \centering
        \includegraphics[width=\textwidth]{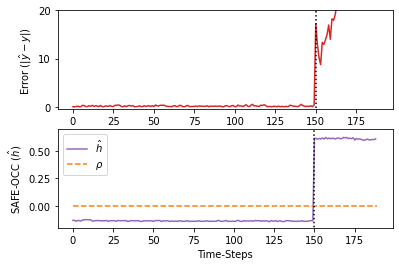}
        \caption{Simulation 3 Novelty Response}
        \label{fig:sim3_error}
    \end{subfigure}
    \quad
    \begin{subfigure}[b]{0.45\textwidth}
        \centering
        \includegraphics[width=\textwidth]{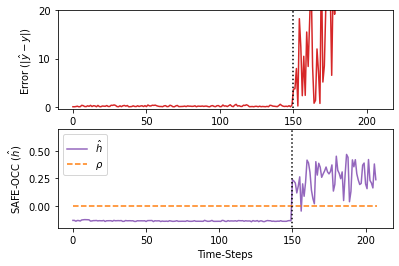}
        \caption{Simulation 4 Novelty Response}
        \label{fig:sim4_error}
    \end{subfigure}
   \caption{Control response trajectories of Simulations 3 and 4. The vertical dotted line at time-step 150 indicates when splattering and blockage disturbances are first introduced to each simulation. The setpoint tracking quickly fails due to the high error incurred by the novel images as corroborated by the SAFE-OCC framework.}
   \label{fig:cart_control_novel}
\end{figure}

Figure \ref{fig:cart_control_novel} shows the response trajectories of Simulations 3 and 4 which are subjected to image splattering and blockage disturbances, respectively. Each of these disturbance types are novel relative to the CNN sensor; thus, significant prediction error is incurred in each case once the sensor is subjected to the disturbance. These highly erred state predictions are injected into the controller which quickly deviates from the set-point until the pendulum completely falls, and the simulation is terminated. In each case, we observe that the SAFE-OCC novelty trajectories accurately identify the novel images once they are injected into the CNN sensor. This highlights how the SAFE-OCC novelty detector effectively identifies novel image data that can incur catastrophic control failure if no recourse action is taken. Hence, in practice the SAFE-OCC novelty detector should be incorporated into a safety system that can take appropriate recourse action once the SAFE-OCC novelty detector identifies novel process data.

\section{Conclusions and Future Work} \label{sec:conclusion}

We have demonstrated that the SAFE-OCC novelty detection framework readily incorporates with CNN sensors to effectively assess the novelty of incoming visual data. This contrasts traditional approaches that conduct novelty detection with independently derived 1D feature spaces. This key difference it what makes the SAFE-OCC framework a natural choice for designing safety systems to mitigate the risk of CNN sensors injecting highly inaccurate measurements into a control system (potentially leading to costly operational deviations). Moreover, SAFE-OCC provide methodological flexibility in how it is implemented and in what machine learning techniques it leverages. This allows us to tailor make novelty detectors in accordance with the unique aspects of a particular process.
\\

In future work, this flexibility should be explored further to better identify the advantages and appropriate use cases of candidate configurations. For instance, the conditions under which certain CNN feature map layers should be selected over others warrants further investigation, as does studying the properties of the candidate scalarization functions that can be applied to the feature map matrices. Furthermore, the extension suggested in Section \ref{sec:feature_occ} of employing an ensemble of SAFE-OCC novelty detectors with varied configurations warrants future investigation, since such an approach could help to greatly increase the sensitivity of the novelty detection and potentially could provide uncertainty quantification. More investigation into effective safety system architectures mitigate the effects of erroneous CNN sensor data once detected by the SAFE-OCC framework also warrants further research. Moreover, applying the SAFE-OCC framework to a real-world control process would be a valuable research direction.

\section*{Acknowledgments}

Most of this work was completed during a summer research internship of Joshua Pulsipher at ExxonMobil.  We acknowledge partial support from the U.S. Department of Energy under grant DE-SC0014114 used to complete the manuscript.

\appendix
\section{Additional Case Study Results} \label{sec:appendix}

\begin{table}
    \centering
    \caption{SAFE-OCC configuration novelty accuracies (\%) for novel inputs to CNN Sensor A.}
    \begin{tabular}{|c|c|c|c|c|c|c|c|c|}
    \hline
        Input & $s(\cdot)$ & $g(\cdot)$ & Splattered & Blurred & Blockages & Fogged & Noised & Shifted \\ \hline \hline
        \multirow{6}{*}{$\Psi^{(1)}$} & \multirow{2}{*}{$-$} & $g_\text{2d2pca}$ & 1.30 & 0.26 & 72.79 & 52.08 & 0.91 & 52.73 \\
         &  & $g_\text{max}$ & 98.44 & 100.00 & 72.53 & 100.00 & 97.79 & 37.37 \\ 
         & \multirow{2}{*}{$s_\text{scale}$} & $g_\text{2d2pca}$ & 100.00 & 0.00 & 94.92 & 100.00 & 100.00 & 58.85 \\ 
         &  & $g_\text{max}$ & 100.00 & 94.01 & 70.96 & 100.00 & 100.00 & 62.11 \\ 
         & \multirow{2}{*}{$s_\text{standard}$} & $g_\text{2d2pca}$ & 100.00 & 6.90 & 96.88 & 100.00 & 100.00 & 89.45 \\ 
         &  & $g_\text{max}$ & 100.00 & 100.00 & 75.65 & 100.00 & 100.00 & 75.00 \\ \hline
        \multirow{6}{*}{$\Psi^{(4)}$} & \multirow{2}{*}{$-$} & $g_\text{2d2pca}$ & 100.00 & 0.52 & 97.01 & 100.00 & 100.00 & 76.04 \\ 
         &  & $g_\text{max}$ & 100.00 & 0.00 & 100.00 & 100.00 & 100.00 & 64.32 \\ 
         & \multirow{2}{*}{$s_\text{scale}$}  & $g_\text{2d2pca}$ & 100.00 & 0.00 & 99.35 & 100.00 & 100.00 & 88.54 \\ 
         &  & $g_\text{max}$ & 100.00 & 0.00 & 100.00 & 100.00 & 100.00 & 80.21 \\ 
         & \multirow{2}{*}{$s_\text{standard}$} & $g_\text{2d2pca}$ & 100.00 & 4.04 & 100.00 & 100.00 & 100.00 & 97.40 \\ 
         &  & $g_\text{max}$ & 100.00 & 0.00 & 100.00 & 100.00 & 100.00 & 90.49 \\ \hline
        \multirow{6}{*}{$P^{(1)}$} & \multirow{2}{*}{$-$} & $g_\text{2d2pca}$ & 100.00 & 18.88 & 87.50 & 96.35 & 72.27 & 92.45 \\
         &  & $g_\text{max}$ & 99.35 & 100.00 & 57.94 & 100.00 & 92.45 & 22.14 \\ 
         & \multirow{2}{*}{$s_\text{scale}$}  & $g_\text{2d2pca}$ & 99.09 & 1.04 & 70.57 & 98.96 & 99.09 & 33.20 \\ 
         &  & $g_\text{max}$ & 100.00 & 97.27 & 70.57 & 100.00 & 100.00 & 61.85 \\ 
         & \multirow{2}{*}{$s_\text{standard}$} & $g_\text{2d2pca}$ & 100.00 & 5.47 & 89.84 & 100.00 & 100.00 & 76.95 \\ 
         &  & $g_\text{max}$ & 100.00 & 100.00 & 74.61 & 100.00 & 100.00 & 74.87 \\ \hline
        \multirow{6}{*}{$P^{(4)}$} & \multirow{2}{*}{$-$} & $g_\text{2d2pca}$ & 100.00 & 0.00 & 98.96 & 100.00 & 100.00 & 72.40 \\ 
         &  & $g_\text{max}$ & 100.00 & 0.00 & 100.00 & 100.00 & 100.00 & 65.10 \\ 
         & \multirow{2}{*}{$s_\text{scale}$}  & $g_\text{2d2pca}$ & 100.00 & 0.00 & 100.00 & 100.00 & 100.00 & 80.34 \\ 
         &  & $g_\text{max}$ & 100.00 & 0.00 & 100.00 & 100.00 & 100.00 & 80.21 \\ 
         & \multirow{2}{*}{$s_\text{standard}$} & $g_\text{2d2pca}$ & 100.00 & 0.00 & 100.00 & 100.00 & 100.00 & 96.35 \\ 
         &  & $g_\text{max}$ & 100.00 & 0.00 & 100.00 & 100.00 & 100.00 & 90.49 \\ \hline
        \multirow{6}{*}{$A^{(1)}$} & \multirow{2}{*}{$-$} & $g_\text{2d2pca}$ & 95.44 & 35.94 & 92.06 & 100.00 & 97.79 & 99.48 \\ 
         &  & $g_\text{max}$ & 99.35 & 100.00 & 57.94 & 100.00 & 92.45 & 22.14 \\ 
         & \multirow{2}{*}{$s_\text{scale}$}  & $g_\text{2d2pca}$ & 6.38 & 0.52 & 65.36 & 59.38 & 0.00 & 29.69 \\ 
         &  & $g_\text{max}$ & 100.00 & 97.27 & 70.57 & 100.00 & 100.00 & 61.85 \\ 
         & \multirow{2}{*}{$s_\text{standard}$} & $g_\text{2d2pca}$ & 84.11 & 3.91 & 83.85 & 100.00 & 98.83 & 64.32 \\ 
         &  & $g_\text{max}$ & 100.00 & 100.00 & 74.61 & 100.00 & 100.00 & 74.87 \\ \hline
        \multirow{6}{*}{$A^{(4)}$} & \multirow{2}{*}{$-$} & $g_\text{2d2pca}$ & 100.00 & 0.00 & 98.70 & 100.00 & 100.00 & 78.52 \\ 
         &  & $g_\text{max}$ & 100.00 & 0.00 & 100.00 & 100.00 & 100.00 & 65.10 \\ 
         & \multirow{2}{*}{$s_\text{scale}$}  & $g_\text{2d2pca}$ & 100.00 & 0.00 & 99.74 & 100.00 & 100.00 & 90.49 \\ 
         &  & $g_\text{max}$ & 100.00 & 0.00 & 100.00 & 100.00 & 100.00 & 80.21 \\ 
         & \multirow{2}{*}{$s_\text{standard}$} & $g_\text{2d2pca}$ & 100.00 & 5.21 & 100.00 & 100.00 & 100.00 & 98.05 \\ 
         &  & $g_\text{max}$ & 100.00 & 0.00 & 100.00 & 100.00 & 100.00 & 90.49 \\ \hline
    \end{tabular}
    \label{tab:CleanImageStudy}
\end{table}

\bibliography{references}
\end{document}